\documentclass[a4paper,10pt]{amsart}
\usepackage{amsmath,amsfonts,amssymb, amsthm, xypic}
\usepackage{bbm}
\usepackage{epsfig,psfrag}
\usepackage[all]{xy}

\def\dim{\mathrm{dim}}
\def\Hom{\mathrm{Hom}}
\def\Z{\mathbb{Z}}
\def\Q{\mathbb{Q}}
\def\QQ{\underline{Coh}}
\def\E{{X}}
\def\mathb{\mathfrak}
\def\Ext{\mathrm{Ext}}
\def\H{\mathbf{H}}
\def\U{\mathbf{U}}
\def\qed{$\hfill \checkmark$}
\def\N{\mathbb{N}}
\def\tto{\twoheadrightarrow}

\def\a{\alpha}
\def\b{\beta}
\def\qlb{\overline{\mathbb{Q}_l}}
\def\ss{\textbf{ss}}
\def\C{\mathbb{C}}
\def\O{\mathcal{O}}

\def\x{\mathbf{x}}
\def\y{\mathbf{y}}

\def\p{\mathbf{p}}
\def\fqb{\overline{\mathbb{F}_q}}
\def\fq{\mathbb{F}_q}
\def\Eb{\overline{X}}
\def\Rb{\mathbf{R}}

\def\LL{\boldsymbol{\Lambda} \hspace{-.06in} \boldsymbol{\Lambda}}
\def\LLambda{\boldsymbol{\Lambda}}
\def\ZZ{\mathbf{Z}}
\def\s{\sigma}
\def\bs{\bar{\sigma}}

\def\UU{\boldsymbol{\mathcal{E}}}
\def\kk{\boldsymbol{{k}}}
\def\kkb{\bar{\boldsymbol{{k}}}}
\def\deg{\textbf{deg}}
\def\rank{\textbf{rank}}
\def\Coh{\textbf{C}}

\newtheorem{theo}{\bf{Theorem}}[section]
\newtheorem{lem}[theo]{Lemma}
\newtheorem{sublem}[theo]{Sublemma}
\newtheorem{cor}[theo]{Corollary}
\newtheorem{prop}[theo]{Proposition}

\numberwithin{equation}{section}
\numberwithin{claim}{section}

\title{On the Hall algebra of an elliptic curve, II}
\author{Olivier Schiffmann}
\thanks{\noindent * Institut Math\'ematique de Jussieu, Paris 6, 175 rue du Chevaleret, 75013 Paris,  FRANCE, 
e-mail:\;\texttt{olive@math.jussieu.fr}}

\begin{document}
\maketitle

{{\flushright{{
\quote{\flushright{\textit{Of you I ask one thing alone,\\
Leave, leave your ancient lore unknown !}\\
H.P. Lovecraft \\
}}}}}}

\tableofcontents
\setlength{\unitlength}{10pt}

\section{Introduction}

\vspace{.1in} 

\paragraph{} This paper is a companion to \cite{BS}, where a categorical approach to the rings of diagonal invariants
$$\LL=\C[x_1^{\pm 1}, \ldots, y_1^{\pm 1}, \ldots]^{\mathfrak{S}_{\infty}}, \qquad
\LL^+=\C[x_1^{\pm 1}, \ldots, y_1, \ldots ]^{\mathfrak{S}_{\infty}}$$
was described in terms of the so-called Hall algebra $\mathbf{U}^+_{\E}$ of the category of coherent sheaves on an elliptic curve $\E$ (defined over a finite field $\mathbb{F}_q$). 
This Hall algebra turns out to be a two-parameter deformation of $\LL^+$, the two deformation parameters being the Frobenius eigenvalues $\s, \bs$ of the particular elliptic curve.  In fact, the structure constants for $\mathbf{U}^+_{\E}$ are Laurent polynomials in $\s$ and $\bs$, and $\mathbf{U}^+_{\E}$ is the specialization of some ``universal'' Hall algebra $\UU^+_{\Rb}$ defined over the ring $\Rb=\C[\s^{\pm 1/2}, \bs^{\pm 1/2}]$. This provides a generalization, to the rings of diagonal invariants, of Steinitz and Hall's well-known realization of the ring of symmetric functions in terms of the classical Hall algebra (see \cite{Mac}, Chap.II).

\vspace{.1in}

The aim of the present work is to turn the above ``categorical'' construction into a purely geometric one. For this, rather than the category $Coh(\E)$ of coherent sheaves on $\E$, we consider the moduli spaces (stacks) $\QQ^{\a}$ of coherent sheaves on $\E$, where $\a$ runs among the set of all possible pairs $(r,d)$ of ranks and degrees. We construct (in a manner reminiscent of \cite{Luszchar}, \cite{L1}, and also of \cite{S1}) a certain category $\mathcal{Q}^{\a}$ of constructible complexes on $\QQ^\a$ along with two convolution functors
\begin{align*}
\text{Ind}^{\a,\beta}&: \mathcal{Q}^{\a} \boxtimes \mathcal{Q}^{\beta} \to \mathcal{Q}^{\a+\beta},\\
\text{Res}^{\a,\beta}&: \mathcal{Q}^{\a+\beta} \to \mathcal{Q}^{\a} \; \hat{\boxtimes}\; \mathcal{Q}^{\beta}.
\end{align*}
This allows us to endow the Grothendieck group $\mathfrak{U}^+_{\E}={\bigoplus}_{\a} K_0(\mathcal{Q}^{\a})$ with an algebra and a coalgebra structure. Moreover, working over $k=\overline{\mathbb{F}_q}$ and using Grothendieck's faisceaux-function correspondence yields a trace morphism
$$Tr_1: \widehat{\mathfrak{U}}^+_{\E} \stackrel{\sim}{\to} \widehat{\mathbf{U}}^+_{\E}$$
which is compatible with the bialgebra structures on both sides. Here $\widehat{\mathfrak{U}}^+_{\E}$, $\widehat{\mathbf{U}}^+_{\E}$ are certain explicit completions of $\mathfrak{U}^+_{\E}$ and ${\mathbf{U}}^+_{\E}$. 

\vspace{.1in}

The collection $\mathcal{P}$ of all simple objects (simple perverse sheaves) in the categories $\mathcal{Q}^{\a}$ can be completely determined and yields, via the above trace map a ``canonical'' basis of $\widehat{\mathbf{U}}^+_{\E}$. In fact, we show that this canonical basis comes from a unique ``universal'' basis $\{{\mathbf{b}}_{\mathbb{P}}\}$ of $\widehat{\UU}^+_{\Rb}$ and that this latter basis admits a purely algebraic characterization in terms of an involution and a lattice. 
This may be viewed as a generalization, in the setting of diagonal invariants, of Lusztig's realization in \cite{LusGreen} of the Schur polynomials as the trace of some simple perverse sheaves on the nilpotent variety of $\mathfrak{gl}(n)$. Motivated by this analogy, we define the \textit{elliptic Kostka polynomials} as the coefficients $\daleth_{\mathbb{P},\mathbb{Q}}(\s,\bs)$ of the transition matrix between the canonical basis $\{{\mathbf{b}}_{\mathbb{P}}\}$ and a natural ``PBW''-type basis $\{\rho_{\mathbb{Q}}\}$ of $\widehat{\UU}^+_{\Rb}$ which plays the role, in our context, of the basis of Hall-Littlewood symmetric functions. By construction, these polynomials satisfy some positivity property, but should not be confused with the $(q,t)$-Kostka polynomials of \cite{Mac}, Chap.VI. The relation between these two families of polynomials is explained in \cite{SV}. We hope and expect that the polynomials $\daleth_{\mathbb{P},\mathbb{Q}}(\s,\bs)$ will play an interesting role in combinatorics and representation theory.

Finally, we note that it is crucial for us to work with the \textit{whole} moduli stack $\QQ^{\a}$ and not just its (semi)stable locus;  indeed it is precisely the singularities of the unstable locus which the polynomials $\daleth_{\mathbb{P},\mathbb{Q}}(\s,\bs)$ describe. 

\vspace{.1in}

To finish this introduction, we state several more possible interpretations of the above constructions. 

\vspace{.05in}
\noindent
First of all, as shown in \cite{SV}, the algebra $\UU^+_{\Rb}$ projects onto the (positive) \textit{spherical double affine Hecke algebra} $\mathbf{S}\ddot{\mathbf{H}}^+_n$ of type $\mathfrak{gl}(n)$ for any $n$. In analogy with the case of the spherical affine Hecke algebras, $\mathbf{S}\ddot{\mathbf{H}}^+_n$ should be isomorphic to a certain convolution algebra of perverse sheaves on a Schubert variety $X_n$ of a ``double affine Grassmanian'' $\ddot{Gr}_n$. Though such an object doesn't exist at the moment, it appears that the stacks $\QQ^{\a}$ provide a model for $X_n$ in the stable limit $n \mapsto \infty$ (to be more precise, the categories $\mathcal{Q}^{\a}$ provide a model for the categories of equivariant perverse sheaves on the limit of $X_n$ as $n \mapsto \infty$). There are however two noteworthy differences between our present situation and the classical picture of the affine Grassmanian~: the presence of simple perverse sheaves associated to nontrivial local systems, and the fact that the Frobenius eigenvalues of the stalks of these perverse sheaves are not all equal to $q^{i/2}$ for some integer $i$ but rather belong to $\s^{\Z}\bs^{\Z}$. This accounts for the fact that $\UU^+_{\Rb}$ depends on two parameters rather than one.

\vspace{.05in}
\noindent
The second interpretation is based on an observation of Loojienga (unpublished, see \cite{GV}).
Let $\mathcal{L}G$ be the \textit{holomorphic} loop group of $GL(n)$ and let ${\mathcal{L}G}\rtimes \C^*$ denotes its one-dimensional universal central extension. Fix $q \in \C^*$ and denote by $\E=\C^*/q^{\Z}$ the associated elliptic curve. Then there is a one-to-one correspondence between conjugacy classes in $(\mathcal{L}G \times \{q\})$ and isomorphism classes of holomorphic vector bundles on $\E$ of rank $n$. Thus, the holomorphic analogs of the categories $\mathcal{Q}^{\a}$ may (heuristically !) be thought of as a category of $\mathcal{L}G$-equivariant perverse sheaves on $(\mathcal{L}G \times \{q\})$. 
For a finite-dimensional reductive group such categories of perverse sheaves were considered and studied in details by Lusztig (the \textit{character sheaves}, see \cite{Luszchar}) and have been shown to be of fundamental importance in representation theory. The $\mathcal{Q}^{\a}$ thus provide a possible model for an extension of Lusztig's character sheaf theory to holomorphic loop groups of type $A$. We thank Victor Ginzburg for kindly explaining to us this point of view.

\vspace{.05in}
\noindent
The final interpretation is in terms of the geometric Langlands program. Recall that this program for the group $GL(n)$ aims at setting up a correspondence between rank $n$ local systems on a smooth projective curve $\E$ defined over a finite field $\mathbb{F}_q$ and a certain collection of perverse sheaves on the moduli stack $\underline{Bun}_n(\E)$ of rank $n$ vector bundles on $X$. In \cite{Lau}, Laumon constructed a category of perverse sheaves on the stacks $\underline{Bun}_n(\E)$ of rank $n$ vector bundles (the so-called  \textit{Eisenstein sheaves}) which should be relevant to the above Langlands correspondence (in the case of the \textit{trivial} local system on $\E$, or more precisely the formal neighborhood of the trivial local system on $\E$). Our categories $\mathcal{Q}^{\a}$ are precisely the categories generated by the simple factors of the Eisenstein sheaves when $\E$ is an elliptic curve. In particular, we obtain a complete description of all the simple factors of the automorphic sheaves in that case, closed formulas for their induction/restriction products, and an algorithm to compute the Poincar\'e polynomial of their cohomology stalks. We refer the reader to \cite{SV2} for more in this direction.

\vspace{.1in}

\paragraph{\textbf{Plan of the paper.}} In Section 1 we recall the main notions and results of \cite{BS} concerning the Hall algebra $\mathbf{U}^+_{\E}$ and its generic version $\UU^+_{\Rb}$.  In Section 2 we give a first--purely algebraic--definition of a canonical basis $\mathbf{B}=\{{\mathbf{b}}_{\mathbf{p}}\}$ of a completion $\widehat{\UU}^+_{\Rb}$ in terms of an involution and a lattice. Section 3 introduces the stacks of coherent sheaves $\QQ^{\a}$ and the convolution functors $\text{Ind}$ and $\text{Res}$. This part closely follows \cite{S1}, which was in turn inspired by \cite{L1}. We also define the category $\mathcal{A}^{\a}$ of semisimple complexes and provide a complete description of the simple perverse sheaves appearing there. In Section 4 we study these simple perverse sheaves, and in particular we prove that they are all pointwise pure. In Section 5 we use the trace map to relate the (completed) Grothendieck group $\widehat{\mathfrak{U}}^+_{\E}={\bigoplus}_{\a} \widehat{K_0}(\mathcal{Q}^{\a})$ and the (completed) Hall algebra
$\widehat{\mathbf{U}}^+_{\E}$. This allows us to define a second canonical basis $\{{\mathbf{b}}_{\mathbb{P}}\}$ of $\widehat{\UU}^+_{\Rb}$ by taking the traces of the simple perverse sheaves in $\bigsqcup_{\a}\mathcal{Q}^{\a}$. We then show the equality of the two canonical bases in Section 6, using some support and degree argument. Finally, in the last section, we give the definition of the elliptic Kostka polynomials $\daleth_{\mathbf{p},\mathbf{q}}(\s,\bs)$ as well as some of their first properties (like $SL(2,\Z)$-invariance).

\vspace{.1in}

\paragraph{\textbf{Notations.}} We consider a smooth elliptic curve $\E$ defined over a finite field $\boldsymbol{{k}}=\mathbb{F}_q$ and we set $\overline{\E}=\E \times_{Spec\; \kk} Spec\; \kkb$. We denote by $Coh(\E)$ and $Coh(\overline{\E})$ the categories of coherent sheaves on $\E$ and $\overline{\E}$ respectively. We fix a line bundle $\mathcal{L}$ of degree one on $\E$ and let $x_0$ be the $\mathbb{F}_q$-rational point of $\E$ such that $\mathcal{L}=\mathcal{O}(x_0)$. We will also use the standard notations for partitions: if $\lambda$ is a partition then $l(\lambda)$ and $|\lambda|$ are the length and size of $\lambda$ respectively. Finally, since we will never consider higher extension groups, we denote $\mathrm{Ext}^1$ simply by $\mathrm{Ext}$.

\pagebreak

\section{Reminder on coherent sheaves over an elliptic curve}

\vspace{.1in}

In this section we very briefly recall various classical results describing the categories $Coh(\E)$ and $Coh(\overline{\E})$.
\vspace{.1in}

\paragraph{\textbf{1.1.}}  Here we let $Y$ stand for either $\E$ or $\overline{\E}$.
Define the slope of a sheaf $\mathcal{F} \in Coh(Y)$ by $\mu(\mathcal{F})=\frac{\textbf{deg}(\mathcal{F})}{\textbf{rank}(\mathcal{F})} \in \mathbb{Q} \cup \{\infty\}$. Recall that a sheaf $\mathcal{F}$ is semistable (resp. stable) if for any $\mathcal{G} \subset \mathcal{F}$ we have
$\mu(\mathcal{G}) \leq \mu(\mathcal{F})$ (resp. $\mu(\mathcal{G}) < \mu(\mathcal{F})$). The category 
$\textbf{C}_{\nu}$ of all semistable sheaves of slope $\nu$ is abelian, artinian, and stable under extensions. We will use the following facts which go back to Atiyah \cite{A}~:

\vspace{.1in}

\paragraph{\textbf{a)}} If $\nu > \mu$ then $\text{Hom}(\textbf{C}_{\nu}, \textbf{C}_{\mu})=\text{Ext}(\textbf{C}_{\mu},\textbf{C}_{\nu})=0$.

\vspace{.1in}

\paragraph{\textbf{b)}} Every sheaf $\mathcal{F}$ possesses a canonical filtration (the Harder-Narasimhan filtration) 
$$0=\mathcal{F}_0 \subset \mathcal{F}_1 \subset \cdots \subset \mathcal{F}_r=\mathcal{F}$$
for which $\mathcal{F}_i/\mathcal{F}_{i-1}$ is semistable of slope, say $\mu_i$, and $\mu_1 > \cdots > \mu_r$.  By the $\text{Ext}$-vanishing property stated above, this filtrations splits. In particular, every indecomposable sheaf belongs to $\textbf{C}_{\nu}$ for some $\nu$.

\vspace{.1in}

\paragraph{\textbf{c)}} For any pair $\nu < \nu'$, let $\Coh_{[\nu,\nu']}$ be the full subcategory whose objects are the sheaves isomorphic to direct sums of indecomposables in $\Coh_{\tau}$ for $\tau \in [\nu,\nu']$. Then $\Coh_{[\nu,\nu']}$ is closed under extensions.

\vspace{.1in}

\paragraph{\textbf{d)}} For any $\nu,\mu \in \mathbb{Q} \cup \{\infty\}$ there is an exact equivalence $\epsilon_{\nu,\mu}: \Coh_{\mu} \stackrel{\sim}{\to} \Coh_{\nu}$.
In particular, any $\Coh_{\mu}$ is equivalent to $\Coh_{\infty}=\mathcal{T}or$, the category of torsion sheaves.
For later use we now give, following \cite{Ku} and \cite{LM}, an explicit construction of such an equivalence $\epsilon_{\mu_2,\mu_1}:~\Coh_{\mu_1} \stackrel{\sim}{\to} \Coh_{\mu_2}$ using the concept of \textit{mutations}. If $\mathcal{F}, \mathcal{G}$ are coherent sheaves we define the left (resp. right) mutation of $\mathcal{G}$ with respect to $\mathcal{F}$ via the following canonical sequences~:
\begin{equation}\label{E:leftmut}
\xymatrix{ \Hom (\mathcal{F}, \mathcal{G}) \otimes \mathcal{F} \ar[r] & \mathcal{G} \ar[r] & L_{\mathcal{F}} \mathcal{G} \ar[r] & 0},
\end{equation}
\begin{equation}\label{E:rightmut}
\xymatrix{ 0 \to \Ext (\mathcal{G}, \mathcal{F})^* \otimes \mathcal{F} \ar[r] & R_{\mathcal{F}}\mathcal{G} \ar[r] & \mathcal{G} \ar[r] & 0}.
\end{equation}
Put
\begin{align*}
\mathcal{F}^{\Vdash}&=\{\mathcal{H} \in Coh(\E)\;|\; \mathrm{Hom}(\mathcal{H},\mathcal{F})=0\;\},\\
\mathcal{F}^{\vdash}&=\{\mathcal{H} \in \mathcal{F}^{\Vdash}\;|\;\ \mathrm{Hom}(\mathcal{F}, \mathcal{H}) \otimes 
\mathcal{F}  \to \mathcal{H} \;\mathrm{is\;a\;monomorphism\;}\}.
\end{align*}
The functor $L_{\mathcal{F}}$ induces
an equivalence of categories 
$\mathcal{F}^{\vdash} \to \mathcal{F}^{\Vdash}$, with inverse given by $R_{\mathcal{F}}$ (see \cite{LM} Theorem~4.4).

\vspace{.1in}

Let $F_m$, $m \in \N$ be the collection of all Farey sequences; that is we have $F_0=\{\frac{0}{1}, \frac{1}{0}\}$ and if $F_n=\{\frac{a_1}{b_1}, \ldots, \frac{a_l}{b_l}\}$ then 
$$F_{n+1}=\{\frac{a_1}{b_1}, \frac{a_1+a_2}{b_1+b_2}, \ldots, \frac{a_i}{b_i}, \frac{a_i+a_{i+1}}{b_i+b_{i+1}} ,\frac{a_{i+1}}{b_{i+1}}, \ldots, \frac{a_l}{b_l}\}.$$
It is a standard fact that every positive rational number belongs to $F_n$ for some $n \gg 0$ and that 
$a_{i+1}b_{i}-a_{i}b_{i+1}=1$ for any $i$ and $n$.
Let us call a stable sheaf $\mathcal{F}$ absolutely simple if $\frac{\deg(\mathcal{F})}{\rank(\mathcal{F})}$ is a reduced fraction (in particular, a torsion sheaf is absolutely simple if and only if it is of degree one).

\begin{prop}[see \cite{LM}] Let $\mu_1=\frac{a}{b}, \mu_2=\frac{c}{d}$ be two consecutive entries of  $F_n$ and let $\mathcal{F} \in \Coh_{\mu_1}$ be an absolutely simple sheaf. Then $R_{\mathcal{F}}$ restricts to an exact equivalence $\Coh_{\mu_2} \stackrel{\sim}{\to} \Coh_{\mu}$ where $\mu=\frac{a+c}{b+d}$. \end{prop}

Using this Proposition we first construct a distinguished absolutely simple object $S_{\mu} \in \Coh_{\mu}$ for all $\mu \in \mathbb{Q}^+ \cup \{\infty\}$ : we put $S_{\infty}=\mathcal{O}_{x_0}, S_0=\mathcal{O}$, and if $S_\mu$ is defined for all entries in $F_n$  and $\{\mu_1,\mu',\mu_2\}$ are consecutive entries in $F_{n+1}$ with $\mu_1,\mu_2 \in F_n$ then we set $S_{\mu'}=R_{S_{\mu_1}} S_{\mu_2}$.  From this we may define an equivalence $\epsilon_{\mu,\infty}: \Coh_{\infty} \stackrel{\sim}{\to} \Coh_{\mu}$ inductively for entries $\mu \in \mathbb{Q}^+$~: $\epsilon_{\infty,\infty}=Id$ and if $\epsilon_{\mu,\infty}$ is known for all $0 \neq \mu \in F_n$ and $\{\mu_1,\mu',\mu_2\}$ are consecutive entries in $F_{n+1}$ with $\mu_1,\mu_2 \in F_n$ then we put $\epsilon_{\mu',\infty}=R_{S_{\mu_1}} \circ 
\epsilon_{\mu_2,\infty}$. Finally, it is easy to check that for $\mu \in \mathbb{Q}^+$ we have
$\epsilon_{\mu,\infty} \simeq ( \cdot \otimes \mathcal{L}^*) \circ \epsilon_{\mu+1,\infty}$; we may thus define unambiguously $\epsilon_{\mu, \infty}= ( \cdot \otimes (\mathcal{L}^*)^{\otimes N}) \circ \epsilon_{\mu+N,\infty}$ for any $\mu \in \mathbb{Q}$. We may also define an inverse equivalence $\epsilon^{-1}_{\mu,\infty}~: \Coh_{\mu} \stackrel{\sim}{\to} \Coh_{\infty}$ in a similar way using left mutations and
we put $\epsilon_{\mu_1,\mu_2}=
\epsilon_{\mu_1,\infty} \circ \epsilon^{-1}_{\mu_2,\infty}$.  We extend the definition of the object
$S_{\mu}$ to an arbitrary $\mu$ by putting $S_{\mu}=(\mathcal{L}^*)^{\otimes N} \otimes S_{\mu+N}$ for $N \gg 0$. Observe that $\epsilon_{\mu_1,\mu_2}(S_{\mu_2})=S_{\mu_1}$ for any $\mu_1, \mu_2$.

\vspace{.1in}

\paragraph{\textbf{e)}} By the class of a sheaf $\mathcal{F}$ we will mean the pair $\overline{\mathcal{F}}=(\rank(\mathcal{F}), \deg(\mathcal{F})) \in \Z^2$. More of the structure of $\mathcal{F}$ is encoded in its HN type, which is defined as $HN(\mathcal{F})=(\overline{\mathcal{H}_1}, \ldots, \overline{\mathcal{H}_r})$, if
$\mathcal{F}=\mathcal{H}_1 \oplus \cdots \oplus \mathcal{H}_r$ where $\mathcal{H}_i$ belongs to $\Coh_{\mu_i}$ and $\mu_1 < \cdots < \mu_r$. We introduce an order on the set of HN types as follows : $((r_1,d_1), \ldots (r_s,d_s)) \preceq ((r'_1,d'_1), \ldots,
(r'_t,d'_t))$ if there exists $l$ such that $(r_{s-i},d_{s-i})=(r'_{t-i},d'_{t-i})$ for $i <l$ while 
$\frac{d_{s-l}}{r_{s-l}} > \frac{d'_{t-l}}{r'_{t-l}}$ or $\frac{d_{s-l}}{r_{s-l}} = \frac{d'_{t-l}}{r'_{t-l}}$ and $d_{s-l} > d'_{t-l}$. For such an order, the maximal HN type of a given class $\alpha \in \Z^2$ is simply $(\a)$ (which corresponds to semistable sheaves of class $\a$).

\vspace{.2in}

\paragraph{\textbf{1.2.}}  We will freely use the definitions and notations of \cite{BS}. We briefly recall the basic notions for the reader's convenience. Let $\H_{\E}$ be the Hall algebra of $\E$ and let $\U^+_\E$ be the spherical subalgebra of $\H_\E$ introduced in \cite{BS}, Section 4. The definition of $\H_{\E}$ and $\U^+_\E$ requires the choice of a square root $v$ of $q$; it is important for us to take $v=-q^{-1/2}$ here.

Recall that $\H_{\E}=\bigoplus_{\mathcal{F} \in Coh(\E)} \C [\mathcal{F}]$ has a basis indexed by isoclasses of objects in $Coh(X)$, and that $\U^+_{\E}$ is the subalgebra generated by elements $\{ \mathbf{1}^{\textbf{ss}}_{\a}\;|\; \a \in \ZZ^+\}$ where
$$\ZZ^+=\{(r,d) \in \Z^2\;|\; r >0\;\text{or}\;r=0, d>0\},$$
$$ \mathbf{1}^{\textbf{ss}}_{\a}=\sum_{\substack{\overline{\mathcal{F}}=\a \\ \mathcal{F} \in \textbf{C}_{\mu(\a)}}} [\mathcal{F}].$$
One also introduces elements $\{ \widetilde{T}_{\a}\;|\; \a \in \ZZ^+\}$ and $\{ \mathbf{1}_{\a}\;|\; \a \in \ZZ^+\}$ satisfying
$$1+\sum_{l \geq 1} \mathbf{1}^{\textbf{ss}}_{l\a_0}s^l=exp\bigg( \sum_{l \geq 1} \widetilde{T}_{l\a_0}s^l\bigg)$$
for any $\a_0=(r,d)$ with $\textbf{g.c.d}(r,d)=1$, and
$$ \mathbf{1}_{\a}=\sum_{\overline{\mathcal{F}}=\a} [\mathcal{F}].$$
Note that $\mathbf{1}_\a$ only belongs to a completion of $\U^+_\E$. The notion of a path $\mathbf{p}=(\x_1, \ldots, \x_l))$ in $\Z^+$ is defined in \cite{BS}, Section 5. The set of all convex paths in $\Z^+$ is denoted $\textbf{Conv}^+$. For any path $\p=(\x_1, \ldots, \x_l))$ we set $\widetilde{T}_\p=\widetilde{T}_{\x_1} \cdots \widetilde{T}_{\x_l}$. Then 
$$\U^+_\E=\bigoplus_{\p \in \textbf{Conv}^+} \C \widetilde{T}_\p.$$

\vspace{.15in}

Put $\Rb=\C[\s^{\pm 1/2}, \bs^{\pm 1/2}]$ where $\s, \bs$ are formal variables. A generic version $\UU^+_\Rb$ of $\U^+_\E$ is defined in \cite{BS}, Section 6, which specializes to $\U^+_\E$ when $\s, \bs$ are set to be equal to the Frobenius eigenvalues in $H^1(\overline{\E}, \qlb)$. The algebra $\UU^+_\Rb$ does not depend on $\E$.  Generic forms $\mathbf{1}^{\textbf{ss}}_{\a}, \tilde{t}_{\a}, \mathbf{1}_{\a}$ of the elements $\mathbf{1}^{\textbf{ss}}_{\a}, \widetilde{T}_{\a}, \mathbf{1}_{\a}$ are also defined for any $\a \in \ZZ^+$, and we have $$\UU^+_\Rb=\bigoplus_{\p \in \textbf{Conv}^+} \Rb \tilde{t}_\p.$$

\vspace{.1in}

Both $\U^+_{\E}$ and $\UU^+_\Rb$ are $\ZZ^+$-graded. The components of degree $\a$ are equal to
$$\U^+_\E[\a]= \bigoplus_{\substack{\p \in \textbf{Conv}^+ \\ wt(\p)=\a}} \widetilde{T}_\p, \qquad \UU^+_\Rb[\a]= \bigoplus_{\substack{\p \in \textbf{Conv}^+ \\ wt(\p)=\a}} \Rb \widetilde{T}_\p$$
where by definition $wt((\x_1, \ldots, \x_l))=\sum \x_i$.

\vspace{.2in}

\section{Algebraic construction of $\mathbf{B}$}

\vspace{.1in}

\paragraph{\textbf{2.1.}} The basis $\{\tilde{t}_{\mathbf{p}}\;|\; \mathbf{p} \in \mathbf{Conv}^+\}$
plays the role of a monomial basis for the algebra $\UU_{\Rb}^+ $.
In analogy with the case of the Hall algebra of a quiver we provide here a tentative definition for a ``canonical basis'' of $\UU^+_{\Rb}$ using an involution and a lattice. We will show in Section~6 that this ``canonical basis''  can also be realized geometrically.

\vspace{.1in}

For our purpose, it is necessary to consider a certain formal completion of $\UU_{\Rb}^+ $. Define an adic valuation $\nu$ on $\UU^+_\Rb$ by $ \nu(\tilde{t}_{\p})=-\mu(\x_1)$ if $\p=(\x_1, \ldots, \x_l) \in \textbf{Conv}^+$. Fix some $C>0$ and denote by $|\;| : u \mapsto C^{- \nu(u)}$ the associated adic norm on $\UU^+_\Rb$. For any $\a \in \ZZ^+$ we let $\widehat{\UU}^+_\Rb[\a]$ be the completion of $\UU^+_\Rb$ with respect to $|\;|$ and set $\widehat{\UU}^+_\Rb=\bigoplus_\a \widehat{\UU}^+_\Rb[\a]$. It is easy to see that, for any fixed $\a$ and $n \in \Z$ there exist  finitely many convex paths $\p=(\x_1, \ldots, \x_l)$ of weight $\a$ satisfying $\mu(\x_1) =\geq n$. Hence there is an identification
$$\widehat{\UU}^+_\Rb[\a] =\prod_{\substack{\p \in \textbf{Conv}^+\\ wt(\p)=\a}} \Rb \tilde{t}_\p.$$
As proved in \cite{BS}, Section 2 the multiplication map is continuous with respect to  $|\;|$ and thus $\widehat{\UU}^+_\Rb$ is an algebra as well.

\vspace{.1in}

Note that $\mathbf{1}_\a \in \widehat{\UU}^+_\Rb$ for any $\a$. For any path $\mathbf{p}=(\x_1, \ldots, \x_l)$ we set ${\mathbf{1}}_{\p}={\mathbf{1}}_{\x_1} \cdots {\mathbf{1}}_{\x_l}$. We also put
${\mathbf{1}}^{\ss}_{\p}={\mathbf{1}}^{\ss}_{\x_1} \cdots {\mathbf{1}}^{\textbf{ss}}_{\x_l}$.
The elements $\{{\mathbf{1}}^{\ss}_{\p}\}_{\p}$ are obtained from $\{\tilde{t}_{\p}\}_{\p}$ by an invertible matrix. Therefore $\{{\mathbf{1}}^{\ss}_{\p}\;|\; \p \in \mathbf{Conv}^+\}$ is also a basis of $\UU^+_{\Rb}$.

\vspace{.2in}

\paragraph{\textbf{2.2.}} We define a weak partial order on $\mathbf{Conv}^+$ as follows. For any $\p=(\x_1, \ldots, \x_r)$ and any slope $\mu \in \mathbb{Q} \cup \{\infty\}$ we put $deg_{\mu}(\p)=\sum_{\mu(\x_i)=\mu} deg(\x_i)$. The symbols $deg_{\geq \mu}(\p), deg_{> \mu}(\p), etc.$ have similar meanings. Now, for two convex paths $\p, \mathbf{q}$ we write $\p \preceq \mathbf{q}$ if there exists $\mu$ such that $deg_{\kappa}(\p)=deg_{\kappa}(\mathbf{q})$ for any $\kappa > \mu$ while $deg_{\mu}(\p)\geq deg_{\mu}(\mathbf{q})$. There is a natural projection map from the set of convex paths $\mathbf{Conv}^+$ to the set of HN types (see Section~1.1.), which simply assigns to a path $\mathbf{p}=(\x_1, \ldots)$ the HN type $HN(\mathbf{p})=(\a_1, \ldots, \a_l)$ with
$$\a_1=\x_1 + \cdots + \x_{i_1}, \; \a_2=\x_{i_1+1} + \cdots + \x_{i_2}, \ldots$$
$$\mu(\x_1)=\cdots=\mu(\x_{i_1}) < \mu(\x_{i_1+1})= \cdots=\mu(\x_{i_2}) < \cdots.$$
The weak order on $\mathbf{Conv}^+$ coincides with the pullback of the order on HN types defined in Section~1.1.e).
We make the following useful observation, which may be deduced from the above together with
1.1.c)~: for any convex paths $\p, \mathbf{q}$ and any slope $\mu$ we have 
\begin{equation}\label{E:observation}
deg_{\geq \mu}(\mathbf{o}) \geq deg_{\geq \mu}(\mathbf{q})
\end{equation}
for any $\tilde{t}_\mathbf{o}$ appearing in the product $\tilde{t}_{\p}\tilde{t}_{\mathbf{q}}$.

Finally, note that we may have $\p \sim \mathbf{q}$ (i.e $\mathbf{p} \preceq \mathbf{q}$ and $\mathbf{q} \preceq \p$) but $\p \neq \mathbf{q}$. In fact for any $\p$ there holds
\begin{equation}\label{E:psimq}
{\mathbf{1}}^{\ss}_{\p} \in \bigoplus_{\mathbf{q} \sim \p}\Rb\tilde{t}_{\mathbf{q}}, \qquad
\tilde{t}_{\p} \in \bigoplus_{\mathbf{q} \sim \p} \Rb{\mathbf{1}}^{\ss}_{\mathbf{q}}.
\end{equation}

\vspace{.1in}

\begin{prop}\label{P:basisone}
The set $\{{\mathbf{1}}_{\p}\;|\; \mathbf{p} \in \mathbf{Conv}^+\}$ is a topological basis
of $\widehat{\UU}^+_{\Rb}$, i.e any element $z \in \widehat{\UU}^+_{\Rb}$ may be written in a unique way as a convergent sum
$z=\sum_{i} a_i {\mathbf{1}}_{\p_i}$ with $a_i \in \Rb$ and $\p_i \in \mathbf{Conv}^+$.
\end{prop}

\noindent
\textit{Proof.} We will prove by induction that there exists a converging sum ${\mathbf{1}}_{\mathbf{q}} = {\mathbf{1}}^{\ss}_{\mathbf{q}} + \sum_{\p \prec \mathbf{q}} a_\p\tilde{t}_{\p}$, with $a_\p \in \Rb$. The statement is obvious if $\mathbf{q}=(\x)$ is of length one. So let us fix some $\mathbf{q}=(\x_1, \ldots \x_r)$ and let us assume that the statement holds for any $\mathbf{q}'=(\x'_1, \ldots, \x'_{r'})$ with $r'<r$. As 
$${\mathbf{1}}_{\x_r}\in{\mathbf{1}}^{\ss}_{\x_r} \oplus \prod_{deg_{>\mu(\x_r)}(\mathbf{o})>0} \Rb\tilde{t}_{\mathbf{o}},$$
we deduce using (\ref{E:observation}) that ${\mathbf{1}}_{\mathbf{q}} \in
{\mathbf{1}}_{(\x_1, \ldots, \x_{r-1})} {\mathbf{1}}^{\ss}_{\x_r} \oplus
\prod_{\mathbf{p} \prec \mathbf{q}}\Rb \tilde{t}_{\p}$.
Using (\ref{E:observation}) again, we have $\tilde{t}_{\p'} {\mathbf{1}}^{\ss}_{\x_r}
\in \prod_{\p \prec \mathbf{q}} \Rb\tilde{t}_{\p}$ for any $\p' \prec (\x_1, \ldots, \x_{r-1})$. By the induction hypothesis it thus follows that 
${\mathbf{1}}_{\mathbf{q}} \in {\mathbf{1}}^{\ss}_{\x_1} \cdots {\mathbf{1}}^{\ss}_{\x_r} \oplus \prod_{\p \prec \mathbf{q}} \Rb\tilde{t}_{\p}$
as desired. \\
Using (\ref{E:psimq}) we have shown that $\{{\mathbf{1}}_{\p}\}$ is related to $\{{\mathbf{1}}^{\ss}_{\p}\}$ by an upper triangular matrix with ones on the diagonal. The Proposition is proved since
$\{{\mathbf{1}}^{\ss}_{\p}\}_{\p}$ is a basis of $\UU^+_{\Rb}$.\qed

\begin{cor}\label{C:unssun}
For any $\mathbf{q} \in \mathbf{Conv}^+$,
$${\mathbf{1}}_{\mathbf{q}} \in {\mathbf{1}}^{\textbf{ss}}_{\mathbf{q}} \oplus \prod_{\p \prec \mathbf{q}} \Rb {\mathbf{1}}^{\textbf{ss}}_{\p},\qquad
{\mathbf{1}}^{\textbf{ss}}_{\mathbf{q}} \in {\mathbf{1}}_{\mathbf{q}} \oplus 
\prod_{\p \prec \mathbf{q}} \Rb{\mathbf{1}}_{\p}.$$
\end{cor}
\noindent
\textit{Proof.} This is a consequence of (\ref{E:psimq}). \qed

\vspace{.2in}

\paragraph{\textbf{2.3.}} By Proposition~\ref{P:basisone}, there exists a unique antilinear involution $z \mapsto \overline{z}$ of $\widehat{\UU}^+_{\Rb}$ such that $\s \mapsto \s^{-1}, \bs\mapsto \bs^{-1}$ and $\overline{{\mathbf{1}}_{\p}}={\mathbf{1}}_{\p}$ for all $\p$. By 
Corollary~\ref{C:unssun}, we have
\begin{equation}\label{E:unssbar}
\overline{{\mathbf{1}}^{\ss}_{\p}} \in {\mathbf{1}}^{\ss}_{\p} \oplus \prod_{\mathbf{q} \prec \mathbf{p}} \Rb{\mathbf{1}}^{\ss}_{\mathbf{q}}.
\end{equation}

\vspace{.1in}

Recall that by \cite{BS} Section 4.2, for any fixed $\y \in \ZZ^+$ satisfying $\deg(\y)=1$, the subalgebra $\UU^{+,(\mu(\y))}_{\Rb}$ generated by $\{\tilde{t}_{r\y}\;|\; r \in \mathbb{N}\}$ is canonically isomorphic to the ring of symmetric polynomials $\LLambda^+ \otimes \Rb= \Rb[x_1, x_2, \ldots]^{\mathfrak{S}_r}$. This isomorphism $i_{\mu(\y)}$ is determined by the condition $i_{\mu(\y)}({\mathbf{1}}^{\ss}_{r\y})=s_r$ with $s_r$ being the Schur function. More generally, if $\lambda=(\lambda_1, \ldots, \lambda_s)$ is any partition we set $\beta_{(\lambda_1 \y, \ldots, \lambda_s \y)}=i_{\mu(\y)}^{-1}(s_{\lambda})$. Finally, for any convex path $\p=(\x_1, \ldots, \x_r)$ with $\mu(\x_1)= \cdots
= \mu(\x_{i_1})< \mu(\x_{i_1+1})= \cdots < \mu(\x_{i_t+1})= \cdots =\mu(\x_r)$ we put
$$\beta_{\p}=\beta_{(\x_{i_1}, \ldots, \x_{i_1})} \cdots \beta_{(\x_{i_t+1}, \ldots, \x_r)}.$$
Observe that $\beta_{\x}={\mathbf{1}}^{\ss}_{\x}$ for any $\x \in \ZZ^+$. Moreover, from
(\ref{E:unssbar}) we see that 
\begin{equation}\label{E:lastsal}
\overline{\beta_{\p}} \in \beta_{\p} \oplus \prod_{\mathbf{q} \prec \p} \Rb \beta_{\mathbf{q}}.
\end{equation}

\vspace{.1in}

We are now at last ready to give the definition of the canonical basis. Let $\Rb^> \subset \Rb$ be the $\C$-linear span of monomials
$\s^a\bs^b$ with $a+b <0$. 

\vspace{.1in}

\addtocounter{theo}{1}

\noindent
\textbf{Definition \thetheo.} For any $\p \in \mathbf{Conv}^+$ we denote by ${\mathbf{b}}_{\p}$ the unique element of $\widehat{\UU}^+_\Rb$ satisfying
$$\overline{{\mathbf{b}}_{\p}}={\mathbf{b}}_{\p}, \qquad {\mathbf{b}}_{\p} \in \beta_{\p} \oplus  \prod_{\p' \prec \p} \Rb^> \beta_{\p'}.$$
We call the set ${\mathbf{B}}=\{{\mathbf{b}}_{\p}\;|\; \p \in \mathbf{Conv}^+\}$ the \textit{canonical basis} of
$\widehat{\UU}^+_{\Rb}$. 

\vspace{.1in}

The existence and uniqueness of ${\mathbf{b}}_{\p}$ results, by a classical argument of Kazhdan and Lusztig, from (\ref{E:lastsal}), together with the fact that all the structure constants in $\UU^+_{\Rb}$ are symmetric in $\s, \bs$. We leave the details to the reader.

\vspace{.1in}

\addtocounter{theo}{1}

\noindent
\textbf{Examples \thetheo.} i) For any $\p=(\x_1, \ldots, \x_r)$ for which $\mu(\x_i)=\infty$ for all $i$ we have
$\mathbf{b}_{\p}=\beta_{\p}$ (indeed, there are no paths $\p'$ of the same weight as $\p$ such that $\p' \prec \p$). In other words, the restriction of the canonical basis to the algebra corresponding to the vertical line in $\ZZ^+$ coincides with the usual canonical basis of $\LLambda^+$ (constructed from the Jordan quiver or the nilpotent variety).\\
ii) For any $\x \in \ZZ^+$ we have ${\mathbf{b}}_{\x}={\mathbf{1}}_{\x}$. Indeed, 
$${\mathbf{1}}_{\x}=\beta_{\x} + \sum_r \sum_{\underset{\mu(\x_1)< \cdots <\mu(\x_r)}{\x_1 + \cdots + \x_r=\x}}\nu^{\sum_{i<j}\langle \x_i,\x_j\rangle} \beta_{(\x_1, \ldots \x_r)}$$
where $\nu=-(\s\bs)^{-1/2}$ and $\langle\;,\;\rangle$ is the Euler form of $\E$.

 \vspace{.2in}

\section{Stacks of coherent sheaves and convolution functors}

\vspace{.1in}

In this section, we employ the method of \cite{S1} to construct a geometric incarnation of $\widehat{\UU}^+_\Rb$ along with ``canonical bases'' which enjoy some integrality and positivity properties. This algebra will turn out to be isomorphic to the spherical Hall algebra $\widehat{\U}^+_{\E}$, and the specialization morphism will map ${\mathbf{B}}$ to the canonical basis constructed in this fashion. This will give us in turn some positivity and integrality properties of ${\mathbf{B}}$.
Whenever possible, we refer to \cite{S1} where a similar construction is given in the context of weighted projective lines of genus one. Until the end of the paper we set $\mathbf{k}=\overline{\mathbb{F}_q}=\overline{\kk}$.

\vspace{.2in}

\paragraph{\textbf{3.1.}} For the notions of algebraic stacks, we refer to the book \cite{LauMo}. We view stacks as sheaves of categories  and work in the fppf topology. Hence to define a stack over a field $l$ it is enough the give the functor of $T$-valued points for any scheme $T$ over $l$. If $\mathcal{C}$ is any category, we write $\langle \mathcal{C} \rangle$ for the category with the same objects but in which the morphisms are the isomorphisms of $\mathcal{C}$. By the descent property of stacks, an algebraic stack is uniquely determined by the corresponding functor of $T$-valued points for affine schemes over $l$.

\vspace{.1in}

For any $\a \in \ZZ^+$, we let $\underline{Coh}^\a$ be the stack of coherent sheaves on $\Eb$ of class $\a$, given by the functor $(Aff/\mathbf{k}) \to Groupoids$
\begin{equation*}
\begin{split}
\texttt{Coh}^{\a}:~T \mapsto &\langle T-{flat,\;coherent\;sheaves\;} \mathcal{F} {\;on\;}T \times \Eb\;{such\;that\;} \\ &\qquad \qquad \qquad \qquad  \overline{\mathcal{F}_{|t}}=\a\;{for\;any\;closed\;point\;} t \in T\rangle
\end{split}
\end{equation*}
The stack $\underline{Coh}^{\a}$ is smooth, locally of finite type, and is an increasing union of smooth open substacks $\QQ_n^{\a}$ defined as follows.
For $\mathcal{E} \in Coh(\Eb)$ and $\a \in \ZZ^+$
consider the functor $\texttt{Quot}_{\mathcal{E}}^{\a}$ from the category of smooth schemes over $\mathbf{k}$ to the category of sets defined by
\begin{equation*}
\begin{split}
\texttt{Quot}_{\mathcal{E}}^{\alpha}~: \Sigma \mapsto&\{\phi:
\mathcal{E} \boxtimes \mathcal{O}_\Sigma
\twoheadrightarrow \mathcal{F}\;|
\mathcal{F}\;is\;a\;coherent\;\Sigma-flat\; sheaf\;on\;\Sigma \times \Eb,\\
& \qquad \qquad\;
\mathcal{F}_{|\sigma}\;is\;
of\;class\;\alpha\;for\;all\;closed\;points\;\sigma \in \Sigma\}
\end{split}
\end{equation*}
In the above, two maps $\phi,\phi'$ are identified if their
kernels coincide. It is a well-known theorem of Grothendieck that $\texttt{Quot}_{\mathcal{E}}^{\alpha}$ is represented by a projective scheme $Quot_{\mathcal{E}}^{\alpha}$ (see e.g. \cite{LP}). In particular, if
$n \in \Z$ and $\alpha \in \ZZ^+$ are such that $\langle [\O(n)], \a\rangle \geq 0$ we set $\mathcal{L}_n^\a=\mathcal{O}(n) \otimes \mathbf{k}^{\langle [\mathcal{O}(n)], \alpha \rangle}$ and put $Quot_{n}^{\a}=Quot_{\mathcal{L}_n^\a}^{\a}$.
The scheme $Quot_{n}^{\a}$ is singular in general; the open subfunctor of $\texttt{Quot}_{\mathcal{L}_n^\a}^{\a}$ defined by
\begin{equation*}
\begin{split}
{}'{\texttt{Quot}}_{\mathcal{L}_n^\a}^{\alpha}~: \Sigma\mapsto &\{(\phi:
\mathcal{L}_n^\a \boxtimes \mathcal{O}_\Sigma
\twoheadrightarrow \mathcal{F}) \in \texttt{Quot}^\a_{\mathcal{L}_n^\a}\;|\; \phi_{*|\sigma}: \mathbf{k}^{\langle [\mathcal{O}(n)], \alpha \rangle} \stackrel{\sim}{\to} \Hom(\mathcal{O}(n), \mathcal{F}_{\sigma})\\
& \qquad \qquad \qquad \qquad \qquad \qquad \qquad \qquad {\;for\;all\;closed\;points}\;\sigma \in \Sigma\}
\end{split}
\end{equation*}
is represented by an open subset $Q_n^\a \subset Quot_{n}^{\a}$. The group $G_n^\a=
\mathrm{Aut}(\mathcal{E}_n^\a) \simeq GL(\langle [\O(n)], \a \rangle)$ naturally acts on $Q_n^\a$.
We will say that a sheaf $\mathcal{F}$ is strictly generated by $\mathcal{O}(n)$ if $\mathcal{F}$ is generated by $\mathcal{O}(n)$ and $\mathcal{F}$ has no direct summand belonging to $\Coh_n$.

\begin{lem}\label{L:1} The scheme $Q_n^\a$ is smooth and the set of $G_n^\a$-orbits is in natural bijection $\mathcal{F} \leftrightarrow \mathbf{O}_{\mathcal{F},n}$ with the set of sheaves $\mathcal{F}$ of class $\a$ strictly generated by $\O(n)$.
\end{lem} 
\begin{proof} See e.g. \cite{LP}, Section~8.2, or \cite{S1}, Section~2.2.
\end{proof}

\vspace{.1in}

We let $\QQ_n^{\a}$ be the quotient stack $Q_n^{\a}/G_n^{\a}$. As $\mathcal{O}(2)$ is generated by $\mathcal{O}$, there are open embeddings of stacks $\QQ_n^{\a} \subset \QQ_{m}^{\a}$ for any $n,m$ with $n \geq m+2$, and $\underline{Coh}^{\a}$ is the limit of the corresponding direct system. For any sheaf $\mathcal{F}$ of class $\a$, we denote by $\mathbf{O}_{\mathcal{F}}=\mathbf{O}_{\mathcal{F},n}/G_n^{\a} $ the locally closed substack of $\QQ^\a$ parametrizing coherent sheaves isomorphic to $\mathcal{F}$ (this definition is independent of $n$ for $n$ sufficiently negative).

\vspace{.15in}

The following remark will be useful.

\begin{lem}\label{L:localsyst} Any $G_n^{\a}$-invariant local system on an orbit $\mathbf{O}_{\mathcal{F},n}$ is constant.
\end{lem}
\noindent
\textit{Proof.} For any $z=(\phi: \mathcal{L}_n^{\a} \tto \mathcal{F}) \in \mathbf{O}_{\mathcal{F},n}$ we have have $Stab_{G_n^{\a}}(z) \simeq \text{Aut}(\mathcal{F})$ (see \cite{S1}, Lemma~2.4.), hence $\mathbf{O}_{\mathcal{F},n} \simeq G_n^{\a} /\text{Aut}(\mathcal{F})$. Thus $G_n^{\a}$-invariant local systems on $\mathbf{O}_{\mathcal{F},n}$ are parametrized by representations of the component group of $\text{Aut}(\mathcal{F})$. We now prove that $\text{Aut}(\mathcal{F})$ is connected.
Let us write $\mathcal{F}=\mathcal{H}_1 \oplus \cdots \oplus \mathcal{H}_r$ with $\mathcal{H}_i$ semistable and $\mu(\mathcal{H}_1) < \cdots < \mu(\mathcal{H}_r$. The group $\text{Aut}(\mathcal{F})$ is an affine fibration over
$\text{Aut}(\mathcal{H}_1) \times \cdots \times \text{Aut}(\mathcal{H}_r)$. Hence it is enough to show that
$\text{Aut}(\mathcal{G})$ is connected for any semistable $\mathcal{G}$. This is easily checked if $\mathcal{G}$ is a torsion sheaf, and follows for an arbitrary $\mathcal{G}$ using the equivalences $\epsilon_{\mu,\infty}$.\qed 

\vspace{.1in}

Let $\mathbb{P}$ be a constructible sheaf of $\qlb$-vector spaces on $\underline{Coh}^{\a}$. By restriction it gives rise, for any $n$, to a $G_n^{\a}$-equivariant constructible sheaf on $Q_n^{\a}$, and hence to a $G_n^{\a}$-equivariant local system $\mathfrak{L}_n$ on $\mathbf{O}_{\mathcal{F},n}$.  By the above Lemma, such a local system is constant. Moreover, by construction, there are canonical maps $\Gamma(\mathbf{O}_{\mathcal{F},m},\mathfrak{L}_m) \to  \Gamma(\mathbf{O}_{\mathcal{F},n},\mathfrak{L}_n)$ for $n \geq m+2$, which are isomorphisms for $n,m \ll 0$. We may now define the \textit{stalk} of $\mathbb{P}$ over $\mathbf{O}_{\mathcal{F}}$ to be
$\mathbb{P}_{|\mathbf{O}_{\mathcal{F}}}=\underset{\longleftarrow}{\text{Lim}}\;\Gamma(\mathbf{O}_{\mathcal{F},n},\mathfrak{L}_n)$ (a $\qlb$-vector space).  

\vspace{.2in}

\paragraph{\textbf{3.2.}} Let ${D}^b(\underline{Z})$ stand for the derived category of constructible $\qlb$-sheaves on an algebraic stack $\underline{Z}$ defined over $\mathbf{k}$. All the stacks $\underline{Z}$ considered in this paper will be increasing unions of open substacks $\underline{Z}_n$, $n \in \Z$, each of which will be a quotient stack $\underline{Z}_n=Q_n/G_n$ of a $\mathbf{k}$-scheme $Q_n$ by a reductive group $G_n$. For such stacks $\underline{Z}$ there exists a six operations formalism, as well as a notion of a dualizing sheaf, and a category of perverse sheaves
(see \cite{LauMo}, Chap. 18.8). In that situation, we let $D^b(\underline{Z})^{ss}$ stand for the category of semisimple $\qlb$-constructible complexes of geometric origin. 
Such a complex $\mathbb{P}$ gives rise by restriction (and is essentially equivalent) to the data of a collection of $G_n$-equivariant semisimple complexes $\mathbb{P}_n \in D^b_{G_n}(Q_n)^{ss}$ for $n \in \Z$, together with certain maps between them satisfying some compatibility conditions. Here $D^b_{G_n}(Q_n)$ is the equivariant derived category of constructible sheaves over $Q_n$, as in \cite{BL}. Given a smooth locally closed substack $\underline{Y} \subset \underline{Z}$ together with a semisimple local system $\mathfrak{L}$ on it we denote by $\mathbf{IC}(\underline{Y}, \mathfrak{L}) \in D^b(\underline{Z})^{ss}$ the associated intersection cohomology sheaf.

\vspace{.2in}

\paragraph{\textbf{3.3.}} Following \cite{L1}, we define functors of induction and restriction on the collection of categories $D^b(\underline{Coh}^{\a})^{ss}$ for $\a \in \ZZ^+$. Consider the diagram
\begin{equation}\label{E:diagind}
\xymatrix{\QQ^\beta \times \QQ^\alpha& \underline{\mathcal{E}}^{\a,\b} \ar[l]_-{p_1} \ar[r]^-{p_2}
&
\QQ^{\alpha +\beta}}
\end{equation}
where the following notations are used~:\\
-$\underline{\mathcal{E}}^{\a,\b}$ is the stack associated to the functor $(Aff/\mathbf{k}) \to Groupoids$ given by
\begin{equation*}
\begin{split}
\texttt{{E}}^{\a,\beta}:~T \mapsto \langle &{exact\;sequences\;} 0 \to \mathcal{G} \to \mathcal{F} \to \mathcal{H} \to 0\;\text{of\;}T-{flat\;coherent\;sheaves}\\
&{\;on\;}T \times \overline{\E} {\;such\;that\;}\overline{\mathcal{G}_{|t}}=\alpha, \overline{\mathcal{F}_{|t}}=\a+\beta\;{for\;any\;closed\;point\;} t \in T\rangle
\end{split}
\end{equation*}
-The 1-morphism $p_1$ is induced by the natural transformation $\texttt{E}^{\a,\beta} \to \texttt{Coh}^{\beta} \times \texttt{Coh}^{\a}$ given on the objects by $ (0 \to \mathcal{G} \to \mathcal{F} \to \mathcal{H} \to 0) \mapsto (\mathcal{H}, \mathcal{G})$,\\
-The 1-morphism $p_2$ is induced by the natural transformation $\texttt{E}^{\a,\beta} \to \texttt{Coh}^{\a+\beta}$ given on the objects by $ (0 \to \mathcal{G} \to \mathcal{F} \to \mathcal{H} \to 0) \mapsto (\mathcal{F})$.

The morphism $p_1$ is smooth with connected fibers (see \cite{S1}, Lemma~3.2), while the morphism $p_2$ is proper (its fibers are isomorphic to certain projective Quot schemes).

\vspace{.15in}

We set 
$$\widetilde{\mathrm{Ind}}^{\b,\a}=
 p_{2!} p_1^*: D^b(\QQ^{\b})^{ss} \boxtimes D^b(\QQ^{\a})^{ss} \to D^b(\QQ^{\a+\b})^{ss}$$
and
$\mathrm{Ind}^{\b,\a}=\widetilde{\mathrm{Ind}}^{\b,\a}[-\langle\beta,\alpha\rangle]$. Recall that the $\underline{Coh}^{\a}$ are locally quotient stacks $Q_n^{\a}/G_n^{\a}$, and all semisimple complexes considered here are assumed to be of geometric origin, so that the equivariant version of the Decomposition
theorem \cite{BBD} (see \cite{BL}) implies that
$\mathrm{Ind}^{\b,\a}$ does indeed take its values in
$D^b(\QQ^{\alpha+\beta})^{ss}$.
The functor $\mathrm{Ind}^{\b,\a}$ is associative i.e. for each triple $\alpha, \beta, \gamma \in \ZZ^+$ there are (canonical) natural transformations 
$\mathrm{Ind}^{\beta+\gamma,\a} \circ \mathrm{Ind}^{\gamma,\b}
\simeq \mathrm{Ind}^{\gamma,\alpha+\beta} \circ
\mathrm{Ind}^{\b,\a}$. This allows us to define an iterated induction functor $\mathrm{Ind}^{\alpha_r, \ldots, \alpha_1}:\; D^b(\QQ^{\alpha_r})^{ss}
\times \cdots \times D^b(\QQ^{\alpha_1})^{ss}
\to D^b(\QQ^{\alpha_1 + \cdots + \alpha_r})^{ss}.$
Note that the functor $\text{Ind}^{\b,\a}$ commutes with Verdier duality.

\vspace{.15in}

Similarly, we set 
$$\widetilde{\text{Res}}^{\beta,\a}=p_{1!}p_{2}^*~:D^b(\QQ^{\a+\b})^{ss} \to D^b({\QQ}^{\b} \times {\QQ}^{\a})$$
and
$\mathrm{Res}^{\b,\a}_n=\widetilde{\mathrm{Res}}_n^{\b,\a}
[-\langle \beta,\alpha\rangle]$. The functor $\text{Res}^{\b,\a}$ in fact restricts to a functor $D^b_{G_n^{\a+\beta}}(Q_n^{\a+\b})^{ss} \to {D}^b_{G_n^{\beta} \times G_n^{\a}}(Q_n^{\b} \times Q_n^{\a})$, whose definition we unravel in detail for later purposes. For this, we fix
a subspace $V \subset \mathbf{k}^{\langle \mathcal{O}(n), \a+\b\rangle}$ of dimension $\langle \mathcal{O}(n),\a \rangle$ along with isomorphisms
$a: V \stackrel{\sim}{\to} \mathbf{k}^{\langle \mathcal{O}(n),\a \rangle}, b:
\mathbf{k}^{\langle \mathcal{O}(n),\a+\b\rangle}/V \stackrel{\sim}{\to}
\mathbf{k}^{\langle \mathcal{O}(n), \b \rangle}$, we put $\underline{V}=V \otimes \mathcal{O}(n)$ and we consider the diagram
\begin{equation}\label{E:12}
\xymatrix{
Q_n^{\alpha+\beta} & F  \ar[l]_-{i} \ar[r]^-{\kappa} & Q_n^\b
\times Q_n^\a}
\end{equation}
where\\
- $F$ is the subvariety of $Q_n^{\alpha+\beta}$ whose points are the quotients $(\phi:
\mathcal{L}^{\alpha+\beta}_n
\tto \mathcal{F})$ such that $\overline{\phi(\underline{V})}=\alpha$
and $i: F \hookrightarrow Q_n^{\alpha+\beta}$ is the (closed) embedding,\\
- $\kappa(\phi)=( b_*
\phi_{|\mathcal{L}^{\alpha+\beta}_n/\underline{V}},a_{*} \phi_{|\underline{V}})$.

By \cite{S1}, Lemma~3.2., $\kappa$ is a vector
bundle of rank $\langle \mathcal{E}^\beta_n-\beta,
\alpha \rangle$.  If $\mathbb{P}$ is an arbitrary complex in $D^b(\QQ^{\a+\b})^{ss}$ and if $\mathbb{P}_n$ denotes its restriction to $\QQ_n^{\a+\b}$ then the restriction of $\widetilde{\text{Res}}^{\b,\a}(\mathbb{P})$ to $\QQ_n^{\b} \times \QQ_n^{\a}$ is isomorphic to $\widetilde{\text{Res}}^{\b,\a}_n(\mathbb{P}_n)$ where by definition
$$\widetilde{\mathrm{Res}}^{\b, \a}_n=\kappa_{!}i^*:
D^b_{G_n^{\a+\b}}(Q_n^{\a+\b})
\to {D}^b_{G_{n^\b \times G_n^\a}}(Q_n^{\b} \times Q_n^\a).$$

As $\mathrm{Res}^{\a,\b}$ does not a priori preserve semisimple complexes, it does not lift to a functor from $D^b(\QQ^{\a+\b})^{ss}$ to $D^b(\QQ^\b)^{ss} \boxtimes D^b(\QQ^\a)^{ss}$. 

\vspace{.2in}

\paragraph{\textbf{3.4.}} The collection of constant complexes $\big(\qlb_{Q_n^\a}[\mathrm{dim}\;Q_n^\a]\big)_{n \in \Z}$ gives rise to a simple perverse sheaf on $\QQ^\a$ which we denote by $\mathbbm{1}_\a$. 
Let $\mathcal{P}^\a$ be the set of all simple objects in $D^b(\QQ^\a)^{ss}$ which appear (up to a shift) in an induction product 
$\mathrm{Ind}^{\a_r,\ldots, \a_1}(\mathbbm{1}_{\a_r} \boxtimes \cdots \boxtimes \mathbbm{1}_{\a_1})$ 
for some $\a_1, \ldots, \a_r$ such that $\sum_i \a_i=\a$ and $rank(\a_i) \leq 1$ for all $i$. We also set
$\mathcal{P}=\bigsqcup_{\a} \mathcal{P}^\a$. We will be concerned here with the full triangulated subcategory $\mathcal{Q}^\a$ of $D^b(\QQ^\a)^{ss}$ whose objects are the complexes isomorphic to a locally finite sum ${\bigoplus}_i \mathbb{P}^i[d_i]$ with $\mathbb{P}^i \in \mathcal{P}^\a$ for all $i$. We also define triangulated categories $\mathcal{Q}^\b \hat{\boxtimes}
\mathcal{Q}^\a$ and $\mathcal{Q}^\b {\boxtimes}
\mathcal{Q}^\a$ whose objects respectively consist of locally finite, resp. finite sums of objects of the form $\mathbb{P}_\b \boxtimes
\mathbb{P}_{\a}$ with
$\mathbb{P}_\b \in \mathcal{Q}^\b$ and $\mathbb{P}_\a \in
\mathcal{Q}^\a$ (see \cite{S1}, Section~4.2).

\vspace{.1in}

\begin{lem} For any $\a, \b \in \ZZ^+$, the induction and restriction functors induce functors
\begin{align*}
\mathrm{Ind}^{\b,\a}&:\; \mathcal{Q}^\b \boxtimes
\mathcal{Q}^\a \to \mathcal{Q}^{\alpha+\beta},.\\
\mathrm{Res}^{\b,\a}&:\; \mathcal{Q}^{\a+\b} \to \mathcal{Q}^\b
\hat{\boxtimes} \mathcal{Q}^\a.
\end{align*}
\end{lem}
\noindent
\textit{Proof.} Identical to \cite{S1}, Lemma~4.1.\qed

\vspace{.1in}

In particular, the collection of categories $\mathcal{Q}^\a, \a \in \ZZ^+$ is stable under the restriction functor. As shown in \cite{S1}, Section~4,  the functor $\text{Res}^{\b,\a}$ is coassociative. We will often need to consider the iterated restriction functor
$$\mathrm{Res}^{\alpha_r,\ldots,\alpha_1}:
\mathcal{Q}^{\alpha_1+\cdots+\alpha_r} \to \mathcal{Q}^{\alpha_r}
\hat{\boxtimes}
\cdots \hat{\boxtimes} \mathcal{Q}^{\alpha_1}.$$

\vspace{.2in}

\paragraph{\textbf{3.5.}} In this section we provide a parametrization and a complete description of the perverse sheaves appearing in $\mathcal{P}$ (that is, we give for each of these simple perverse sheaves a corresponding smooth locally closed subvariety along with an irreducible local system on it). 

\vspace{.15in}

\paragraph{\textbf{3.5.1}} We first introduce certain stratifications of the stacks
$\underline{Coh}^{\a}$. Recall from Section~1.1. e). that the HN type of a sheaf $\mathcal{F}$ with HN filtration $0 \subset \mathcal{F}_1 \subset \cdots \subset \mathcal{F}_{r-1} \subset \mathcal{F}$ is  $HN(\mathcal{F})=(\overline{\mathcal{F}_1}, \overline{\mathcal{F}_2/\mathcal{F}_1}, \ldots \overline{\mathcal{F}/\mathcal{F}_{r-1}}) \in (\ZZ^+)^r$. The stack $\QQ^\a$ admits a locally finite stratification by locally closed substacks
$$\QQ^\a = \bigsqcup_{(\a_1, \ldots, \a_r)} \underline{HN}^{-1}(\a_1, \ldots, \a_r)$$
where $\underline{HN}^{-1}(\a_1, \ldots, \a_r) \subset \QQ^\a$ is the substack parametrizing sheaves
$\mathcal{F} \in Coh(\E)$ whose HN type is $(\a_1, \ldots, \a_r)$.
For simplicity, the open substack $\underline{HN}^{-1}(\a) \subset \QQ^\a$ corresponding to semistable sheaves will be denoted $\QQ^{(\a)}$. Note
that $\underline{Coh}^{\a} \setminus \QQ^{(\alpha)}=\bigcup_{\underline{\beta}
\prec \alpha} \underline{HN}^{-1}(\underline{\beta})$, where $\prec$ stands for the order defined in Section~1.1.e).

\vspace{.1in}

Any $\a \in \ZZ^+$ may be written in a unique way as $\a=l \delta_{\mu(\a)}$ where $l \geq 1$ and $\delta_{\mu(\a)}=(p,q) \in \ZZ^+$ with $p,q$ relatively prime (so that $\mu(\a)=\frac{p}{q}$). If $\mu \in \Q \cup \{\infty\}$
 then $\QQ^{(\delta_{\mu})}$ actually corresponds to stable sheaves and $\QQ_n^\a$ is a (smooth) geometric quotient $Q_n^{(\delta_{\mu})}  / G_n^{\delta_{\mu}}$ for $n \ll 0$. Arguing (using mutations) in the same way as in \cite{S1}, Section~10., we obtain, for each $\mu_1, \mu_2 \in \mathbb{Q} \cup \{\infty\}$ a canonical isomorphism $\rho_{\mu_1,\mu_2}:  \QQ^{(\delta_{\mu_2})}  \stackrel{\sim}{\to}  \QQ^{(\delta_{\mu_1})}$. In particular there is an isomorphism $\rho_{\infty,\mu}:~\QQ^{(\delta_{\mu})} \stackrel{\sim}{\to}
 \QQ^{((0,1))} \simeq \Eb/\mathbf{k}^*$, where the multiplicative group $\mathbf{k}^*$ acts trivially. In a similar vein, fix $l \geq 1$ and $\mu \in \Q \cup \{\infty\}$ and let us consider the open substack $\underline{U}^{(l\delta_{\mu})}$ of $\QQ^{(l\delta_{\mu})}$ parametrizing semistable sheaves $\mathcal{F}$ isomorphic to direct sums of stable sheaves in $\Coh_{\mu}$ with distinct support~:
 $\mathcal{F} \simeq \epsilon_{\mu,\infty}(\bigoplus_{i=1}^l \mathcal{O}_{x_i}),\; x_i \in \Eb, x_i \neq x_j \; \text{if}\; i \neq j.$
The stack $\underline{U}^{(l\delta_{\mu})}$ may also be realized as a geometric quotient $Y_n^{(l\delta_{\mu})}/G_n^{l\delta_{\mu}}$ for some smooth open subscheme $Y_n^{(l\delta_{\mu})}$ of $Q_n^{l\delta_{\mu}}$ (for $n \ll 0)$ and there is a canonical isomorphism $\rho_{\infty,\mu}^l: \underline{U}^{(l\delta_\mu)} \stackrel{\sim}{\to} (S^l\Eb \setminus \underline{\Delta}) / (\mathbf{k}^*)^l$ where $\underline{\Delta}=\{(x_i)\;|
x_i=x_j\; \text{for\;some}\; i \neq j\}$ and again $(\mathbf{k}^*)^l$ acts trivially. As a consequence, there is a projection $\pi_1(\underline{U}^{(l\delta_\mu)}) \tto \mathfrak{S}_l$ where $\mathfrak{S}_l$ is the symmetric group on $l$ letters respectively. 

\vspace{.1in}

Finally, if $\mathbb{P}$ belongs to $\mathcal{P}^\a$ then there exists a unique HN type $(\a_1, \ldots \a_r)$ such that $supp(\mathbb{P}) \subset 
\overline{\underline{HN}^{-1}(\a_1, \ldots , \a_r)}$ and $supp(\mathbb{P}) \cap \underline{HN}^{-1}(\a_1, \ldots, \a_r)$ is nonempty; we write $HN(\mathbb{P})=(\a_1, \ldots, \a_r)$ and call $(\a_1, \ldots, \a_r)$ the \textit{generic HN type} of $\mathbb{P}$. The set of all generically semistable $\mathbb{P}$ of weight $\a$ will be denoted $\mathcal{P}^{(\a)}$.

\vspace{.15in}

\paragraph{\textbf{3.5.2.}} Any representation $\sigma$ of $\mathfrak{S}_l$ gives rise to a local system on $S^l \Eb \setminus \underline{\Delta}$ and hence to a local system $\mathfrak{L}_{\sigma}$ on $\underline{U}^{(l\delta_{\mu})}$. The following proposition is proved in the same fashion as Proposition~9.7 in \cite{S1}.

\begin{prop}\label{P:descript-perv} For any $\mu \in \mathbb{Q} \cup \{\infty\}$ and any $l \geq 1$ we have
$$\mathcal{P}^{(l\delta_{\mu})}=\{\mathbf{IC}(\underline{U}^{(l\delta_{\mu})},\sigma)\;|\;\sigma \in Irrep\;\mathfrak{S}_l\}.$$
Furthermore, for each $\a$ there is a canonical bijection
$$\theta_{\a}:\;\mathcal{P}^\a \stackrel{\sim}{\to} \bigsqcup_{\underset{\mu(\a_1)< \cdots < \mu(\a_r)}{\a_1 + \cdots + \a_r=\a}} \mathcal{P}^{(\a_1)} \times \cdots \times \mathcal{P}^{(\a_r)}$$
such that if $\theta_{\a}(\mathbb{P})=(\mathbb{P}^1, \ldots, \mathbb{P}^r)$ then $\mathrm{Ind}(\mathbb{P}^1 \boxtimes \cdots \boxtimes \mathbb{P}^r)\simeq \mathbb{P} \oplus \mathbb{P}'$ with $supp(\mathbb{P}') \subset supp(\mathbb{P})$ and $\text{dim\;}supp(\mathbb{P}') < \text{dim\;}supp(\mathbb{P})$. 
Finally, every $\mathbb{P} \in \mathcal{P}$ is self-dual, i.e $D(\mathbb{P})=\mathbb{P}$ where $D$ is the Verdier duality functor.
\end{prop}

\vspace{.2in}

\section{Purity}

\vspace{.1in}

\paragraph{\textbf{4.1.}} Recall that $\kk=\fq$, $\mathbf{k}=\fqb$ and that $\Eb$ is equal to $\E \times_{\small{Spec(\kk)}}Spec(\mathbf{k})$ for some smooth elliptic curve $\E$ defined over $\kk$. The functor ${}^0\texttt{Quot}_{\mathcal{L}_n^{\a}}^{\a}$ from the category of smooth $\fq$-schemes to sets defined in the same manner as ${\texttt{Quot}}_{\mathcal{L}_n^{\a}}^{\a}$ but replacing $\Eb$ by $\E$, is represented by a projective $\fq$-scheme ${}^0Quot_{\mathcal{L}_n^{\a}}^{\a}$ and there is a canonical identification $Quot_{\mathcal{L}_n^{\a}}^{\a}\simeq {}^0Quot_{\mathcal{L}_n^{\a}}^{\a}\times_{\small{Spec(\kk)}}Spec(\mathbf{k})$. A similar statement holds for the open subschemes ${}^0Q_n^{\a} \subset {}^0Quot_{\mathcal{L}_n^\a}^{\a}$ and $Q_n^\a \subset Quot_{\mathcal{L}_{n}^\a}^{\a}$.  In particular, there is a natural (geometric) Frobenius automorphism
$\tilde{F}: Q_n^\a \to Q_n^\a$.  As a consequence, there are quotient stacks ${}^0\QQ_n^{\a}$ defined over $\kk$ such that $\QQ_n^{\a} \simeq {}^0\QQ_n^{\a} \times_{\small{Spec(\kk)}}Spec(\mathbf{k})$. These stacks form an open cover of a stack ${}^0\QQ^{\a}$ defined over $\kk$ and we have $\underline{Coh}^{\a} \simeq {}^0\QQ^{\a} \times_{\small{Spec(\kk)}}Spec(\mathbf{k})$.  This equips $\QQ^\a$ with a geometric Frobenius automorphism still denoted $\tilde{F}$.

\vspace{.1in}

We refer to \cite{Del} and \cite{FW} for matters concerning the notions of Weil sheaves, of purity and of pointwise purity for complexes of constructible $\qlb$-sheaves on a scheme. We will say that a complex $\mathbb{P} \in D^b(\QQ^{\a})^{ss}$ equipped with a mixed structure $h: \mathbb{P} \stackrel{\sim}{\to} \tilde{F}^* \mathbb{P}$ is pure (resp. pointwise pure) of weight $l$ if for any $n \in \Z$ the corresponding $G_n^{\a}$-equivariant complex
$\mathbb{P}_n$ over $Q_n^{\a}$ is pure (resp. pointwise pure) of weight $l$. 
Note that for any coherent sheaf $\mathcal{F} \in Coh(\E)$ the action of the Frobenius on the stalks of $\mathbb{P}_n$ over $\mathbf{O}_{\mathcal{F},n}$ for all $n$ gives rise to a similar Frobenius action on the stalk
$\mathbb{P}_{|\mathbf{O}_{\mathcal{F}}}$ (see Section~3.1.).
A mixed complex $\mathbb{P} \in D^b(\QQ^{\a})^{ss}$ is pointwise pure of weight zero if and only if for any $\mathcal{F} \in Coh(\E)$ the eigenvalues of $\tilde{F}^*$ on the stalks $H^i(\mathbb{P})_{|\mathbf{O}_{\mathcal{F}}}$ are all of
complex norm $q^{i/2}$. 
 We let $\qlb(1/2)$ stand for the square root of the Tate sheaf, whose Frobenius map has eigenvalue $q^{-1/2}$. We write $\mathbb{P}(n/2)$ for $\mathbb{P} \otimes (\qlb(1/2))^{\otimes n}$. From now on the definitions of the functors $\text{Ind}$ and $\text{Res}$ are understood to include the Tate twist as well, i.e we replace everywhere the shifts $[n]$ appearing in Section~3 by $[n](n/2)$.

\vspace{.1in}

Our aim in this section is to establish the following purity statement.

\begin{theo}\label{P:purity} For any $\a \in \ZZ^+$ and any $\mathbb{P} \in \mathcal{P}^\a$ there exists a (unique up to a scalar) isomorphism $h_{\mathbb{P}}: \mathbb{P} \stackrel{\sim}{\to} \tilde{F}^* \mathbb{P}$, with respect to which $\mathbb{P}$ is pointwise pure of weight zero. 
\end{theo}

\noindent
\textit{Proof.} By definition, the complexes $\mathbbm{1}_{\a}$ are equipped with a mixed structure which is pointwise pure of weight zero. 
From the definition of the induction product and \cite{Del}, Prop. 6.2.6, $\mathrm{Ind}^{\a_1, \ldots, \a_r}(\mathbbm{1}_{\a_1} \boxtimes \cdots \boxtimes \mathbbm{1}_{\a_r})$ is also (globally) pure of weight zero for any $\a_1, \ldots, \a_r$.
The key to Theorem~\ref{P:purity} is to show the next statement. 

\begin{prop}\label{L:purityun} For any $\a_1, \ldots, \a_r$ the complex $\mathrm{Ind}^{\a_1, \ldots, \a_r}(\mathbbm{1}_{\a_1} \boxtimes \cdots \boxtimes \mathbbm{1}_{\a_r})$ is pointwise pure of weight zero. 
\end{prop}
\noindent
\textit{Proof.}  To ease the notation let us put $\mathbbm{1}_{\a_1, \ldots, \a_r}=\text{Ind}^{\a_1, \ldots, \a_r}(\mathbbm{1}_{\a_1} \boxtimes \cdots \boxtimes \mathbbm{1}_{\a_r})$.
Set $\a=\sum \a_i$ and let $\mathcal{F} \in Coh(\Eb)$ be a sheaf of class $\a$. We will work over
each open substack $Q_m^{\a}/G_m^{\a}$. There exists $e$ such that for any $z=(\phi: \mathcal{L}_m^{\a} \tto \mathcal{F}) \in \mathbf{O}_{\mathcal{F},m}$ we have $(\tilde{F}^*)^e(z)=z$. We need to show the following property~:

\vspace{.05in}

\noindent
(a) \textit{For any $z \in \mathbf{O}_{\mathcal{F},m}$, the eigenvalues of} $(\tilde{F}^*)^e$ \textit{on the stalk} $H^i(\mathbbm{1}_{\a_1,\ldots,\a_r})_{|z}$ \textit{are all of complex norm} $q^{ei/2}$.  

\vspace{.05in}

Unraveling the definitions, one has that the stalk at $z$ of the complex $\mathbbm{1}_{\a_1, \ldots, \a_r}$ is  equal to the stalk at $z$ of $p_{!}(\qlb_{E''})$ where~: $E''$ is the variety of pairs $(\phi, V_r \subset \cdots \subset V_1 =\mathbf{k}^{\langle \mathcal{O}(m), \a\rangle})$ such that $(\phi: \mathcal{L}_m^{\a} \tto \mathcal{F}) \in \mathbf{O}_{\mathcal{F},m}$, $V_i/V_{i+1}$ is of dimension $\langle \mathcal{O}(m), \alpha_i \rangle$, and $\phi(V_i \otimes \mathcal{O}(m))/\phi(V_{i+1} \otimes \mathcal{O}(m))$ is a coherent sheaf of class $\alpha_i$; $p$ is the projection on the first factor. Hence the fiber of $p$ at $z$ is identified with the projective (hyperquot) scheme $Quot_{\mathcal{F}}^{\a_2, \ldots, \a_r}$ parametrizing successive quotients $\psi: \mathcal{F}= \mathcal{G}_1 \tto \cdots \tto \mathcal{G}_{r}$ with $\mathcal{G}_i$ of class $\a_i + \cdots + \a_r$ (see \cite{S1}, Lemma~4.2). In particular, the stalk of $H^i( p_{!}(\qlb_{E''}))$ at $z$ is isomorphic to $H_c^i(Quot_{\mathcal{F}}^{\a_2, \ldots, \a_r}, \qlb)$ and statement (a) for a sheaf $\mathcal{F}$ is equivalent to

\vspace{.05in}

\noindent
(a') \textit{The hyperquot scheme $Quot_{\mathcal{F}}^{\a_2, \ldots, \a_r}$ is cohomologically pure, i.e the eigenvalues of $(\tilde{F}^*)^e$ on $H^i_c({Quot}_{\mathcal{F}}^{\a_2, \ldots, \a_r})$ are all of complex norm} $q^{ei/2}$.  

\vspace{.05in}

Note that by base change it is enough to prove this for any $\mathcal{F}$ for which $e=1$ and we will assume this from now on. We will prove (a) (or (a')) in three steps. First, we reduce (a) for an arbitrary sheaf $\mathcal{F}$ to (a) for its semistable subquotients, then to (a) for stable sheaves, and finally we prove (a) for stable sheaves.
In the course of the proof we will often need the following lemma~:

\begin{lem}\label{L:s151} Let $\a_1, \ldots, \a_r, \beta, \sigma \in \ZZ^+$ be such that $\sum \a_i=\beta+\sigma$. Then
$$\text{Res}^{\sigma,\beta}(\mathbbm{1}_{\a_1, \ldots, \a_r})=\bigoplus_{\underline{\beta}, \underline{\sigma}} \mathbbm{1}_{\sigma_1, \ldots, \sigma_r} \boxtimes \mathbbm{1}_{\beta_1, \ldots, \beta_r}[u(\underline{\sigma},\underline{\beta})](u(\underline{\sigma},\underline{\beta})/2)$$
where $(\underline{\beta}=(\beta_1, \ldots, \beta_r), \underline{\sigma}=(\sigma_1, \ldots, \sigma_r))$ run through the set of all tuples satisfying $\a_i=\sigma_i+\beta_i, \sum \beta_i=\beta, \sum\sigma_i=\sigma$, and where $u(\underline{\sigma},\underline{\beta})$ are certain integers.\end{lem}

\noindent
\textit{Proof}. Identical to \cite{S1}, Lemma~4.1. (see also \cite{S1}, Corollary~4.3).

\vspace{.1in}

\textit{Step 1.} Let $\mathcal{F}$ be an arbitrary coherent sheaf on $\Eb$ of class $\a$ and let us write $\mathcal{F}=\mathcal{H}_1 \oplus \cdots \oplus \mathcal{H}_s$ where $\mathcal{H}_i$ are semistable sheaves such that $\mu(\mathcal{H}_1)\leq \mu(\mathcal{H}_2) \cdots \leq \mu(\mathcal{H}_s)$. We may and will further assume that for any $i$, $\epsilon_{\infty,\mu(\mathcal{H}_i)}(\mathcal{H}_i)$ is a torsion sheaf supported at a single point $x_i \in \Eb$ and that $x_i \neq x_j$ if $\mu(\mathcal{H}_i) = \mu(\mathcal{H}_j)$ but $i \neq j$. Now let $\mathbb{P} \in D^b(\QQ^\a)^{ss}$ be any complex and let us consider the stalk of $\mathrm{Res}^{\overline{\mathcal{H}_1}, \ldots, \overline{\mathcal{H}_s}}(\mathbb{P})$ over a point $\underline{x}=(x_1, \ldots, x_s)$ with $x_i \in \mathbf{O}_{\mathcal{H}_i,m} $. Setting $\overline{\mathcal{H}_i}=\beta_i$ and using the notation in the restriction diagram
\begin{equation}\label{E:Proppur1}
\xymatrix{
Q_m^{\alpha} & F  \ar[l]_-{i} \ar[r]^-{\kappa} & Q_m^{\b_1} \times \cdots
\times Q_m^{\b_s}}
\end{equation}
we have $\text{Res}^{\b_1, \ldots, \b_s}(\mathbb{P})_{|\underline{x}}=j_{\underline{x}}^*\kappa_{!}i^*(\mathbb{P})=\kappa_{!}(j')^*(\mathbb{P})$ where $j_{\underline{x}}: \{\underline{x}\} \to Q_m^{\b_1} \times \cdots \times Q_m^{\b_r}$ and
$j': \kappa^{-1}(\underline{x}) \to Q_m^{\a}$ are the embeddings. By construction, $\text{Ext}(\mathcal{H}_i, \mathcal{H}_j)=0$ if $i<j$ and hence $\kappa^{-1}(\underline{x}) \subset \mathbf{O}_{\mathcal{F},m}$. As $\mathbb{P}_m$ is $G_m^{\a}$-equivariant, its restriction to $\mathbf{O}_{\mathcal{F},m}$ is constant by Lemma~\ref{L:localsyst}. As $\kappa$ is a vector bundle, we deduce from this that $\text{Res}^{\b_1, \ldots, \b_r}(\mathbb{P})_{|\underline{x}}\simeq\mathbb{P}_{|z}[-2rk(\kappa)](-rk(\kappa))$ for any $z \in \mathbf{O}_{\mathcal{F},m}$. Hence $\mathbb{P}$ is pure at the point $z$ if and only if $\text{Res}^{\b_1, \ldots, \b_s}(\mathbb{P})$ is pure at the point $\underline{x}$. In particular, taking $\mathbb{P}=\mathbbm{1}_{\a_1, \ldots, \a_r}$ we obtain that $\mathbbm{1}_{\a_1, \ldots, \a_r}$ is pure at any $z \in \mathbf{O}_{\mathcal{F},m}$ if and only if
$\mathrm{Res}^{\b_1, \ldots, \b_s}(\mathbbm{1}_{\a_1, \ldots, \a_r})$ is pure at $\underline{x}=(x_1, \ldots , x_s)$. But by Lemma~\ref{L:s151} the stalk at $\underline{x}$ of $\text{Res}^{\b_1, \ldots, \b_s}(\mathbbm{1}_{\a_1, \ldots, \a_r})$ is a sum of complexes which are each, up to shift, an external product of stalks of the form $(\mathbbm{1}_{\delta_1, \ldots, \delta_r})_{|x_i}$. Therefore (a) for $\mathcal{F}$ is a consequence of (a) for each of the semistable sheaves $\mathcal{H}_i$.

\vspace{.1in}

\textit{Step 2.} Let $\mathcal{F}$ be a semistable sheaf of slope $\mu$ and assume that $\epsilon_{\infty,\mu}(\mathcal{F})$ is a torsion sheaf supported at a single point $x \in \Eb$. Thus $\mathcal{F}$ is an iterated extension of a stable sheaf $\mathcal{H}$. Set $\delta_{\mu}=[\mathcal{H}]$ so that $[\mathcal{F}]=l\delta_{\mu}$ for some $l \in \N$. Assume that (a) holds for the stable sheaf $\mathcal{H}$. Then arguing as in Step 1 we see that $\text{Res}^{\delta_{\mu}, \ldots, \delta_{\mu}}(\mathbbm{1}_{\a_1, \ldots, \a_r})$ is pure over any point $\underline{x}=(x_1, \ldots, x_l)$ with $x_i \in \mathbf{O}_{\mathcal{H},m}$. In order to deduce from this that $\mathbbm{1}_{\a_1, \ldots, \a_r}$ is itself pure over any point $z \in \mathbf{O}_{\mathcal{F},m}$, we consider for any $n >0$ the closed subscheme $Q_{m,x}^{n\delta_{\mu}} \subset Q_{m}^{(n\delta_{\mu})}$ parametrizing quotients
$(\psi:\mathcal{L}_m^{n\delta_{\mu}} \tto \mathcal{G})$ for which $\mathcal{G}$ is semistable and $\epsilon_{\infty,\mu}(\mathcal{G})$ is supported at $x$. The quotient stack $Q_{m,x}^{n\delta_{\mu}}/G^{n\delta_{\mu}}_m$ is isomorphic to the quotient stack $\mathcal{N}_n/G_n$ where $\mathcal{N}_n \subset \mathfrak{gl}(n,k)$ is the nilpotent cone and $G_n=GL(n,k)$ acts by conjugation. The restriction functor $\text{Res}^{\delta_{\mu}, \ldots, \delta_{\mu}}$ induces a functor
$T: D^b_{G_l}(\mathcal{N}_l) \to D^b_{G_1 \times \cdots \times G_1}(\mathcal{N}_1 \times \cdots \times \mathcal{N}_1)$.

\begin{lem} Let $\mathbb{P} \in D^b_{G_l}(\mathcal{N}_l)$ be a semisimple complex equipped with a mixed structure. Assume that
$T(\mathbb{P})$ is pointwise pure of weight zero. Then $\mathbb{P}$ is also pointwise pure of weight zero.\end{lem}
\noindent
\textit{Proof.} Since any direct summand of a pointwise pure complex is again pointwise pure, it is enough to consider the case of a simple complex $\mathbb{P}=\mathbf{IC}(\mathcal{O}_{\lambda})[c](d/2)$, where $\mathcal{O}_{\lambda}$ is a nilpotent orbit and $c, d$ are integers. By \cite{KazhLusztig} Theorem~5.5 (see also \cite{LusGreen}), $\textbf{IC}(\mathcal{O}_{\lambda})$ is pointwise pure of weight zero. Moreover, we have $T(\mathbf{IC}(\mathcal{O}_{\lambda}))\simeq (\qlb \boxtimes \cdots \boxtimes \qlb)^{\oplus d_{\lambda}}$ where $d_{\lambda}$ is the multiplicity of the irreducible $gl(l)$-module $V_{\lambda}$ of highest weight $\lambda$ in the tensor product $V_{(1)} \otimes \cdots \otimes V_{(1)}$. Hence $T(\mathbb{P})$ is pointwise pure if and only if $c=d$. But then $\mathbb{P}$ is itself pointwise pure. \qed

\vspace{.05in}

Now, the restriction of $\mathbbm{1}_{\a_1, \ldots, \a_r}$ to $Q_{m,x}^{l\delta_{\mu}}$ is a semisimple complex and since we have assumed that (a) holds for stable sheaves, it satisfies the conditions of the above Lemma. Hence it is pointwise pure as desired.

\vspace{.1in}

\textit{Step 3.} Finally, we deal with the case of a stable sheaf $\mathcal{F}$. If $\mathcal{F}$ is a stable torsion sheaf then $\mathcal{F}=\mathcal{O}_x$ for some $x \in \Eb$ and $Quot_{\mathcal{F}}^{\a_2, \ldots, \a_r}$ is either empty or reduced to a point so that (a') clearly holds. 
Let us now suppose that $\mathcal{F}$ is a line bundle.  In that case, the stalk 
$(\mathbbm{1}_{\a_1, \ldots, \a_r})_{|z}$ is nonzero only if $\a_2, \ldots, \a_{r}$ are torsion classes.  For any line bundle $\mathcal{L}$ and any torsion class $\beta$, the scheme $Quot_{\mathcal{L}}^{\b}$ is isomorphic to the symmetric product $S^{deg(\beta)}\Eb$. It follows that $Quot_{\mathcal{F}}^{\a_2, \ldots, \a_r}$ is an iterated fibre bundle with fibres $S^{deg(\a_i)}\Eb$. Thus it is smooth projective and by \cite{DeligneWeilI}, Theorem I.6, its cohomology is pure and (a') holds. We now argue by induction on the rank of $\mathcal{F}$. Let $\mathcal{F}$ be a stable sheaf of rank $r>1$ and assume that the stalk of any complex $\mathbbm{1}_{\b_1, \ldots, \b_s}$ is pure over a point $(\phi: \mathcal{L}_m^{\a} \tto \mathcal{G})$ whenever $\mathcal{G}$ is a stable sheaf of rank less than $r$. By Atiyah's construction (see Section 1.1 d)) $\mathcal{F}$ is the universal extension of two stable sheaves $\mathcal{G}, \mathcal{H}$ of smaller rank satisfying $\text{dim}\;\text{Ext}(\mathcal{G},\mathcal{H})=1$. Consider the complex $\mathbb{R}=\text{Res}^{\overline{\mathcal{G}},\overline{\mathcal{H}}}(\mathbbm{1}_{\a_1, \ldots, \a_r})$. By Lemma~\ref{L:s151} and the induction hypothesis the stalk $\mathbb{R}_{|(x,y)}$ is pure of weight zero when $x \in \mathbf{O}_{\mathcal{G},m}$ and $y \in \mathbf{O}_{\mathcal{H},m}$.  By definition, $\mathbb{R}=\kappa_{!}i^*(\mathbbm{1}_{\a_1, \ldots, \a_r})$ where
\begin{equation}
\xymatrix{
Q_m^{\alpha} & F  \ar[l]_-{i} \ar[r]^-{\kappa} & Q_m^{\overline{\mathcal{G}}} 
\times X_m^{\overline{\mathcal{H}}}}
\end{equation}
is the restriction diagram. Since $\mathcal{F}$ is the only nontrivial extension of $\mathcal{G}$ by $\mathcal{H}$ we have $\kappa^{-1}((x,y)) \subset \mathbf{O}_{\mathcal{F},m} \cup \mathbf{O}_{\mathcal{G} \oplus \mathcal{H},m}$. Moreover $V:=\kappa^{-1}(x,y)$ is a vector space of dimension $d:=\langle \overline{\mathcal{L}^{{\mathcal{G}}}_m}-\overline{\mathcal{G}},\overline{\mathcal{H}}\rangle$ and $W:=\kappa^{-1}(x,y) \cap \mathbf{O}_{\mathcal{G} \oplus \mathcal{H},m}$ is a vector subspace of $\kappa^{-1}(x,y)$ of dimension $\langle \overline{\mathcal{L}^{{\mathcal{G}}}_m},\overline{\mathcal{H}}\rangle=d-1$. Thus from the decomposition $V=W \sqcup (V \setminus W)$ we deduce a long exact sequence in cohomology
with compact support
$$
\xymatrix{
\cdots \ar[r] &H_c^{l-1}(j_{!}j^*\mathbbm{1}_{\a_1, \ldots, \a_r}) \ar[r]^-{x_{l-1}} & H^l_c(\iota_{!}\iota^*\mathbbm{1}_{\a_1, \ldots, \a_r}) \ar[r]^-{y_l} & H^l(\mathbb{R})_{|(x,y)} \ar[r] & \cdots}$$
where $\iota: V\setminus W \to V$ and $j: W \to V$ are the embeddings. Put $\mathbb{P}=\iota^* \mathbbm{1}_{\a_1, \ldots, \a_r}, \mathbb{P}'=j^*\mathbbm{1}_{\a_1, \ldots, \a_r}$. By Lemma~\ref{L:localsyst}, $\mathbb{P}$ and $\mathbb{P}'$ are constant so that 
$$H^{l-1}_c(j_{!}j^*\mathbbm{1}_{\a_1, \ldots, \a_r})=\bigoplus_h H^{l-1-h}_c(W) \otimes H^h(\mathbb{P}')\simeq H^{l+1-2d}(\mathbb{P}')(2d-2)$$
and 
\begin{equation*}
\begin{split}
H^l_c(\iota_{!}\iota^*\mathbbm{1}_{\a_1, \ldots, \a_r})&=\bigoplus_h H^{l-h}_c(V\setminus W) \otimes H^h(\mathbb{P})\\
&\simeq H^{l-2d}(\mathbb{P})(2d) \oplus H^{l-2d+1}(\mathbb{P})(2d-2).
\end{split}
\end{equation*}
In addition since we have assumed that (a) holds for $\mathcal{G}, \mathcal{H}$ it follows from Step 1. that (a) holds also for $\mathcal{G} \oplus \mathcal{H}$ and hence $\mathbb{P}'$ is pure.  Now, since
$\mathbbm{1}_{\a_1, \ldots, \a_r}$ is globally pure of weight zero, the Frobenius weights in $H^{l-2d+1}(\mathbb{P})(2d-2)$ are all at most $l-1$. On the other hand $\mathbb{R}_{(x,y)}$ is pure hence the weights in $H^l(\mathbb{R})_{(x,y)}$ are all equal to $l$ and therefore $H^{l-2d+1}(\mathbb{P})(2d-2) \subset \text{Ker}\;y_l=\text{Im}\; x_{l-1}$. But all the weights in
$H^{l+1-2d}(\mathbb{P}')(2d-2)$ are equal to $l-1$ as $\mathbb{P}'$ is pure and hence the weights
in $H^{l-2d+1}(\mathbb{P})(2d-2)$ are all equal to $l-1$ so that $\mathbb{P}$ is itself pure. By definition this means that $\mathbbm{1}_{\a_1, \ldots, \a_r}$ is pure over $\mathbf{O}_{\mathcal{F},m}$ and the induction argument is closed. Thus Step 3. is finished, and Proposition~\ref{L:purityun} is proved.\qed

\vspace{.2in}

\noindent
\textit{Proof of Theorem~\ref{P:purity}.}
Let us first assume that $\mathbb{P} \in \mathcal{P}^{(\a)}$. By Proposition~\ref{P:descript-perv}, $\mathbb{P}_m$ is of the form $\mathbf{IC}(Y_m^{(\a)}, \mathfrak{L}_\sigma)$ for some  irreducible $\mathfrak{S}_l$-module $\sigma$. Since $\mathfrak{L}_{\sigma}$ is the extension to $\mathbf{k}$ of a local system ${}^0\mathfrak{L}_{\sigma}$ on the $\kk$-scheme ${}^0Y^{(\a)}_m$, it follows from \cite{FW} Section III, Cor. 9.2 that the unique isomorphism
$h_{\mathbb{P}_m}: \mathbb{P}_m \stackrel{\sim}{\to} \tilde{F}^*\mathbb{P}_m$ whose restriction to the stalk of an $\kk$-point 
of ${}^0Y^{(\a)}_m$ is the identity endows $\mathbb{P}_m$ with a mixed structure which
is globally pure of weight zero.  In the case of a general $\mathbb{P} \in \mathcal{P}^\a$ we may, by Proposition~\ref{P:descript-perv}, identify
$\mathbb{P}$ with an isotypical component of a product $\mathrm{Ind}(\mathbb{P}^1 \boxtimes \cdots \boxtimes \mathbb{P}^r)$ for certain simple perverse sheaves $\mathbb{P}^i \in \mathcal{P}^{(\a_i)}$. By the above, each $\mathbb{P}^i$ is equipped with a mixed structure, globally pure of weight zero. From the definition of the induction product and \cite{Del}, Prop. 6.2.6 it follows that $\mathrm{Ind}(\mathbb{P}^1 \boxtimes \cdots \boxtimes \mathbb{P}^r)$ is also pure of weight zero. The same holds for its isotypical components, and hence for $\mathbb{P}$.

\vspace{.05in}

We now turn to the pointwise purity property. By construction, each $\mathbb{P} \in \mathcal{P}$ appears in some induction product $\mathrm{Ind}^{\a_1, \ldots, \a_r}(\mathbbm{1}_{\a_1} \boxtimes 
\cdots \boxtimes \mathbbm{1}_{\a_r})$. By Proposition~\ref{L:purityun} above, the complex $\mathbb{R}=\mathrm{Ind}^{\a_1, \ldots, \a_r}(\mathbbm{1}_{\a_1}\boxtimes \cdots \boxtimes \mathbbm{1}_{\a_r})$ is pointwise pure of weight zero. On the other hand, it decomposes as $\mathbb{R}=
\bigoplus_{\mathbb{P}} \mathbb{V}_{\mathbb{P}} \otimes \mathbb{P}$, where $\mathbb{V}_{\mathbb{P}}=
 \bigoplus_j\mathrm{Hom}(\mathbb{P},\;^p \hspace{-.03in}H^j(\mathbb{R}))[-j]$ and $\;^p \hspace{-.03in}H^j$ denotes perverse cohomology. Restricting to the stalk over a point $x \in Q_m^\a$ for some $m \in \Z$, we deduce that for any $i$,
$H^i(\mathbb{R}_m)_{|x}=\bigoplus_{h+l=i}\bigoplus_{\mathbb{P}} H^h(\mathbb{V}_{\mathbb{P}}) \otimes H^l(\mathbb{P}_m)_{|x}$ (as $\tilde{F}$-modules).
 Since $\mathbb{P}$ and $\mathbb{R}$ are both globally pure, the complex $\mathbb{V}_{\mathbb{P}}$, viewed as a sheaf over a point, is pure of weight zero. The purity of $(\mathbb{P}_m)_{|x}$ is now a consequence of the purity of $\mathbb{V}_{\mathbb{P}}$ and of $(\mathbb{R}_m)_{|x}$. 
The Theorem is proved. \qed

\vspace{.2in}

\section{The trace map}

\vspace{.1in}

\paragraph{\textbf{5.1.}} We may at this point use Theorem~\ref{P:purity} to construct a collection of (topological) bialgebras $\widehat{\mathfrak{U}}^{+}_{\E,e}$ over $\C$ indexed by the set of positive integers $e \in \mathbb{N}$. Let us fix an isomorphism $ \qlb \simeq \C$, and for any complex of finite-dimensional $\qlb$-vector spaces $V^\bullet$ equipped with an action of $\tilde{F}$ let us set $tr_e(V)=\sum_i (-1)^itr_{H^i(V)}(\tilde{F}^e) \in \C$. 
We consider the $\C$-vector space $\mathfrak{U}^{+}_{\E,e}:= \bigoplus_{\mathbb{P} \in \mathcal{P}} \C \mathb{b}_{\mathbb{P}}$ with a basis $\{\mathb{b}_{\mathbb{P}}\}_{\mathbb{P}}$ indexed by $\mathcal{P}$. We say that an element
$\mathb{b}_{\mathbb{P}}$ is of $h$-degree $n$ if
$\mathbb{P}=(\mathbb{P}_m)_{m \in \Z}$ with $\mathbb{P}_m=0$ for $m > -n$ and
$\mathbb{P}_{-n} \neq 0$. Denote by
$\widehat{\mathfrak{U}}_{\E,e}^{+}$
the completion of
$\mathfrak{U}^{+}_{\E,e}$ with respect to the $h$-adic topology.
There is a natural multiplication map  $\mathfrak{U}^{+}_{\E,e}
\otimes\mathfrak{U}^{+}_{\E,e} \to
\widehat{\mathfrak{U}}^{+}_{\E,e}$ defined by
$\mathb{b}_{\mathbb{P}'} \mathb{b}_{\mathbb{P}''}=\sum_{\mathbb{P} \in \mathcal{P}} tr_e(\mathbb{V}_{\mathbb{P}}) \mathb{b}_{\mathbb{P}}$ where $\mathrm{Ind}(\mathbb{P}' \boxtimes \mathbb{P}'')=\bigoplus_{\mathbb{P}} \mathbb{V}_{\mathbb{P}} \otimes \mathbb{P}$. Here, the multiplicity spaces are equipped with a mixed structure coming from that on $\mathbb{P}, \mathbb{P}', \mathbb{P}''$ and 
on $\mathrm{Ind}(\mathbb{P}' \boxtimes \mathbb{P}'')$.
In a similar way, there is a comultiplication map 
$\mathfrak{U}^{+}_{\E,e} \to
\widehat{\mathfrak{U}}^{+}_{\E,e} \hat{\otimes}
\widehat{\mathfrak{U}}^{+}_{\E,e}$ defined by $\Delta(\mathb{b}_{\mathbb{P}})=\sum_{\mathbb{P}', \mathbb{P}''} tr_e(\mathbb{V}^{\mathbb{P}',\mathbb{P}''}) \mathb{b}_{\mathbb{P}'} \otimes \mathb{b}_{\mathbb{P}''}$, where
$\mathrm{Res}(\mathbb{P}) =\bigoplus_{\mathbb{P}', \mathbb{P}''} \mathbb{V}^{\mathbb{P}',\mathbb{P}''}\otimes( \mathbb{P}' \boxtimes \mathbb{P}'')$. Moreover, these maps are continuous (see \cite{S1}, Section~10.1.) and we may extend them
to 
$$m_e: \widehat{\mathfrak{U}}^{+}_{\E,e}
\otimes\widehat{\mathfrak{U}}^{+}_{\E,e} \to
\widehat{\mathfrak{U}}^{+}_{\E,e},$$
$$\Delta_{e}: \widehat{\mathfrak{U}}^{+}_{\E,e} \to
\widehat{\mathfrak{U}}^{+}_{\E,e} \hat{\otimes}
\widehat{\mathfrak{U}}^{+}_{\E,e}.$$

The associativity and coassociativity of $m_e$ and $\Delta_e$ follow from the analogous properties of the functors $\mathrm{Ind}$ and $\mathrm{Res}$ (see \cite{S1}, Section~3).

\vspace{.1in}
 
\paragraph{\textbf{Remark.}} The choice of the identification $\qlb \simeq \C$ is not essential~: as all complexes considered here are pure, the eigenvalues of Frobenius are all algebraic integers and in addition the set of eigenvalues is invariant under $Gal(\overline{\mathbb{Q}}/\mathbb{Q})$. In particular, the algebra $\widehat{\mathfrak{U}}_{\E,e}^{+}$ is independent of that choice.

\vspace{.1in}

For any fixed elliptic curve $\E$ we define the completions $\widehat{\mathbf{H}}_{\E}$ and $\widehat{\mathbf{U}}^+_{\E}$ of $\mathbf{H}_{\E}$ and $\mathbf{U}^+_{\E}$ respectively in the same way as this was done for $\UU^+_{\Rb}$ in Section~2.1 (these completions are also used in \cite{BS}, Section~2.).
Recall that for any coherent sheaf $\mathcal{F} \in Coh(\overline{\E})$ of class $\a$ there is a substack $\mathbf{O}_{\mathcal{F}} \subset \underline{Coh}^{\a}$ parametrizing sheaves isomorphic to $\mathcal{F}$, and that to any complex $\mathbb{P} \in D^b(\QQ^{\a})^{ss}$ is associated its stalk $\mathbb{P}_{|\mathbf{O}_{\mathcal{F}}}$, which is a complex of $\qlb$-vector spaces (see Section~3.1). The same holds for any sheaf $\mathcal{F} \in Coh(\E)$ if we replace $\underline{Coh}^{\a}$ by ${}^0\QQ^{\a}$. Moreover, if $\mathcal{F} \in Coh(\E)$ then $\mathbf{O}_{\mathcal{F}} \subset \underline{Coh}^{\a}$ is (pointwise) fixed by $\tilde{F}$ and if $\mathbb{P} \in \mathcal{P}^{\a}$ is equipped with the mixed structure provided by Theorem~\ref{P:purity} then there is an action of $\tilde{F}$ on the cohomology stalks $H^i(\mathbb{P}_{|\mathbf{O}_{\mathcal{F}}})$. This allows us to define a $\C$-linear trace function $Tr_1: \widehat{\mathfrak{U}}^+_{\E,1} \to \widehat{\mathbf{H}}_{\E}$ given by
$$Tr_1(\mathb{b}_{\mathbb{P}})=\sum_{\mathcal{F}} tr_1(H^{\bullet}(\mathbb{P}_{|\mathbf{O}_{\mathcal{F}}}))[\mathcal{F}].$$
By definition, $Tr_1(\mathb{b}_{\mathbb{P}})$ is equal to the limit as $n$ tends to $-\infty$ of the elements
$$Tr^{(n)}_1(\mathb{b}_{\mathbb{P}})=v^{-\dim\;G_n^\a}\sum_{\mathbf{O} \in Q_n^\a/G_n^\a} tr_1(H^\bullet(\mathbb{P}_{n\;|\mathbf{O}})) [\mathcal{F}_{\mathbf{O}}]
$$
where $v=-q^{-1/2}$ (the factor $v^{-\dim\;G_n^\a}$ guarantees that for any $m < n$, $Tr^{(n)}_1(\mathb{b}_{\mathbb{P}})$ and
$Tr^{(m)}_1(\mathb{b}_{\mathbb{P}})$ coincide on the set of objects of $Coh(\E)$ strictly generated by $\mathcal{O}(n)$).  By Grothendieck's trace formula, this map is a bialgebra homomorphism (see e.g. \cite{SLecturesII}, Section~3.6.). Replacing $1$ by $e \in \N$ throughout yields a similar homomorphism $Tr_e:\widehat{\mathfrak{U}}^+_{\E,e} \to \widehat{\mathbf{H}}_{\E/\mathbb{F}_{q^e}}$, where $\E/\mathbb{F}_{q^e}=\E \times_{Spec\;\kk} Spec\; \mathbb{F}_{q^e}$. Note that in defining $\widehat{\mathbf{H}}_{\E/\mathbb{F}_{q^e}}$ and $Tr_e$ we use the parameters $q^e$ and $v_e=-{q^{-e/2}}$ instead of $q,v$.

\vspace{.1in}

For any $\a \in \ZZ^+$ let us set $\mathb{b}_{\a}=\mathb{b}_{\mathbbm{1}_{\a}} \in \widehat{\mathfrak{U}}^+_{\E,e}$. By construction, we have $Tr_e(\mathb{b}_{\a})=\mathbf{1}_{\a}$.
In the same vein, there exists a unique collection of elements $\{\mathb{b}^{ss}_{\a}\;|\a \in \ZZ^+\}$ 
of $\widehat{\mathfrak{U}}^+_{\E,e}$ satisfying 
\begin{equation}\label{E:ssslip}
\mathb{b}_{\a}=\mathb{b}_{\a}^{ss} + \sum_{t > 1} \sum_{\underset{\mu(\a_1)< \cdots < \mu(\a_t)}{\a_1 + \cdots + \a_t=\a}} v_e^{\sum_{i<j}\langle \a_i,\a_j\rangle}\mathb{b}_{\a_1}^{ss} \cdots
\mathb{b}_{\a_t}^{ss},
\end{equation}
for any $\a$ and we have $Tr_e(\mathb{b}^{ss}_{\a})=\mathbf{1}^{ss}_{\a}$.

\vspace{.05in}

The next Theorem is proved in Section~5.2. 

\begin{theo}\label{T:cano} The subalgebra of $\widehat{\mathfrak{U}}^{+}_{\E,e}$ generated by $\mathb{b}_{\a}$ for $\a\in \ZZ^+$ is dense (in the $h$-adic topology).\end{theo}

\begin{cor}\label{C:cano} The map $Tr_e$ restricts to an isomorphism $\widehat{\mathfrak{U}}^{+}_{\E,e} \stackrel{\sim}{\to} \widehat{\mathbf{U}}^+_{\E/\mathbb{F}_{q^e}}$.\end{cor}
\noindent
\textit{Proof.} Since $Tr_e$ is continuous and $Tr_e(\mathb{b}_{\a}) \in \widehat{\mathbf{U}}^+_{\E/\mathbb{F}_{q^e}}$, the image of $Tr_e$ belongs to 
$\widehat{\mathbf{U}}^+_{\E/\mathbb{F}_{q^e}}$ by Theorem~\ref{T:cano}, and in fact 
$Tr_e(\widehat{\mathfrak{U}}^{+}_{\E,e})=\widehat{\mathbf{U}}^+_{\E/\mathbb{F}_{q^e}}$ by Proposition~\ref{P:basisone}. It remains to show the injectivity of $Tr_e$. Observe that if $\a=(d,0)$ is a torsion class then by Proposition~\ref{P:descript-perv} and the definition of $\widehat{\mathbf{U}}^+_{\E/\mathbb{F}_{q^e}}$, we have $\dim\; \widehat{\mathfrak{U}}^{+}_{\E,e} [\a]=p(d)=\dim\;\widehat{\mathbf{U}}^+_{\E/\mathbb{F}_{q^e}}[\a]$, where $p(d)$ stands for the number of partitions of $d$. Thus the restriction of $Tr_e$ to
$\bigoplus_d \widehat{\mathfrak{U}}^{+}_{\E,e} [(d,0)]$ is injective. More generally, 
let $\mathfrak{U}^{+,\geq n}_{\E,e} \subset \widehat{\mathfrak{U}}^{+}_{\E,e}$ be the subspace linearly spanned by elements $\mathb{b}_{\mathbb{P}}$ such that $\mathbb{P}=(\mathbb{P}_l)_l$ with $\mathbb{P}_n \neq 0$, and
put $\mathbf{U}^{+,\geq n}_{\E/\mathbb{F}_{q^e}}=\U^{+}_{\E/\mathbb{F}_{q^e}}/( \U^{+}_{\E/\mathbb{F}_{q^e}} \cap
\mathbf{H}^{\not\geq n}_{\E/\mathbb{F}_{q^e}})$ where
$$\mathbf{H}^{\not\geq n }_{\E/\mathbb{F}_{q^e}}=\bigoplus_{\mathcal{F} \in \mathcal{I}_{\not\geq n}} \C [\mathcal{F}], $$
$$ \mathcal{I}_{\not\geq n}=\{\mathcal{F}\;|\mathcal{F} \text{\;is\;not\;strictly\;generated\;by\;} \mathcal{O}(n)\}.$$
Composition with $Tr_e$ gives rise to a map $\widehat{\mathfrak{U}}^{+,\geq n}_{\E,e} \to \mathbf{U}^{+,\geq n}_{\E/\mathbb{F}_{q^e}}$, which is still surjective. On the other hand, by Proposition~\ref{P:descript-perv}, we have
$$\dim\; \widehat{\mathfrak{U}}^{+,\geq n}_{\E,e}[\a]=\sum_{\underset{\mu(\a_1)>\cdots > \mu(\a_r) > n}{\a_1+\a_2+\cdots+\a_r=\a}} p(deg(\a_1)) \cdots p(deg(\a_r))$$
and from \cite{BS}, Theorem~5.4. and Lemma~5.6, one sees that $\dim\;\mathbf{U}^{+,\geq n}_{\E/\mathbb{F}_{q^e}}[\a]$ is given by the same formula. Hence the restriction of $Tr_e$ to
$\widehat{\mathfrak{U}}^{+,\geq n}_{\E,e}$ is injective. Passing to the limit as $n$ tends to $-\infty$, we obtain that
$Tr_e$ is injective.\qed

\vspace{.2in}

\paragraph{\textbf{5.2.}} In this Section we give the proof of Theorem~\ref{T:cano}.  It is quite similar to the proof of Theorem~4.ii) in \cite{S1}, but we provide the details for the reader's convenience.
Recall that for $\nu \in \mathbb{Q} \cup \{\infty\}$ we put $\delta_{\nu}=(p,q) \in \ZZ^+$ where $deg((p,q))=1$ and $p/q=\nu$. We extend by bilinearity the notation $\mathb{b}_{\mathbb{P}}$ to an arbitrary (semisimple) complex $\mathbb{S}$ belonging to some $\mathcal{Q}^\a$.

\vspace{.1in}

By \cite{S1}, Proposition~5.1. and Corollary~5.2., we have
\begin{equation}\label{E:semismall}
\begin{split}
{\mathrm{Ind}}^{n_1\delta_{\infty}, \cdots
,n_r\delta_{\infty}}(\mathbf{IC}&(\underline{U}^{(n_1\delta_{\infty})},{\sigma_1})
\boxtimes \cdots \boxtimes
\mathbf{IC}(\underline{U}^{(n_r\delta_{\infty})},{\sigma_r}))\\
&=\mathbf{IC}(\underline{U}^{(l\delta_{\infty})},ind_{\mathfrak{S}_{n_1}\times
\cdots \times
\mathfrak{S}_{n_r}}^{\mathfrak{S}_l}(\sigma_1 \otimes \cdots \otimes
\sigma_r))
\end{split}
\end{equation}
where $l=\sum n_i$. By construction, $\mathb{b}_{l \delta_{\infty}}=\mathbf{IC}(\underline{U}^{(l\delta_\infty)},Id) \in \widehat{\mathfrak{U}}^{+}_{\E,e}$ for all $l \geq 1$. It is well-known that the Grothendieck group of $\mathfrak{S}_l$ is linearly spanned (over $\mathbb{Z}$) by the class of the trivial representation $Id$ and the classes of the induced representations $ind_{\mathfrak{S}_{n_1} \times \cdots \times \mathfrak{S}_{n_r}}^{\mathfrak{S}_l} (\sigma_1 \otimes \cdots\otimes \sigma_r)$ for all tuples $(n_i)_i$ such that $\sum_i n_i=l$ and $n_i <l$ for all $i$. It thus follows by induction from (\ref{E:semismall}) that $\mathb{b}_{\mathbf{IC}(\underline{U}^{(l\delta_{\infty}}),\sigma)} \in \widehat{\mathfrak{U}}^{+}_{\E,e}$ for any $\sigma$. This proves Theorem~\ref{T:cano} for the weights of the form $l \delta_{\infty}$.

\vspace{.1in}

 Let $\mathfrak{W}$ denote the subalgebra of $\mathfrak{U}^{+}_{\E,e}$ generated by $\{\mathb{b}_{\a}\;|\; \a \in \ZZ^+\}$. We also let $\widehat{\mathfrak{U}}^+_{<n}$ be the completion of ${\bigoplus}_{\mathbb{R},\; \mathbb{R}_{n}=0} \C \mathb{b}_{\mathbb{R}}$. Note that for any $\a$ and $n$, the space $\widehat{\mathfrak{U}}_{\E,e}^+[\a]/\widehat{\mathfrak{U}}^+_{<n}[\a]$ is finite-dimensional.
To prove Theorem~\ref{T:cano} for a general weight $\a$, we have to show the following.

\vspace{.05in}

\noindent
\textit{a)\;For any} $\alpha \in \ZZ^+$, $n \in \Z$ \textit{and for any}
$\mathbb{P} \in \mathcal{P}^{\alpha}$ for which $\mathbb{P}_n \neq 0$
\textit{there exists} $u \in \mathfrak{W}$ such that
$\mathb{b}_{\mathbb{P}} \equiv u\;(\text{mod}\; \widehat{\mathfrak{U}}^+_{<n})$.

\vspace{.05in}

We will prove $a)$ by induction on the rank of $\a$.  The case of $\a$ of rank zero is treated above, so let us choose $\mathbb{P} \in \mathcal{P}^\a$ with $rank(\a) \geq 1$ and let us assume that $a)$ is proved for all $\b$ with $rank(\b) < rank(\a)$. Next, we fix $n \in \Z$, and argue by induction on the degree of $\a$, and finally on the HN type of $\mathbb{P}$. Hence we further assume that $a)$ is proved for all $\mathbb{P}'$ of weight $\b$ with $rank(\b)=rank(\a)$ and $deg(\b) < deg(\a)$ and for all $\mathbb{P}' \in \mathcal{P}^\a$ such that $HN(\mathbb{P}') \prec HN(\mathbb{P})$. Note that for any given $n$ and $\mathbb{P}$ there exist only finitely many such $\mathbb{P}' $ as above for which $\mathbb{P}'_n \neq 0$ so that an induction argument is indeed justified.

\vspace{.1in}

Let us write $HN(\mathbb{P})=\underline{\a}=(\a_1, \ldots, \a_r)$, and let us first suppose that $\mathbb{P}$ is not generically semistable, i.e that $r >1$.
Consider the complex
$$\mathbb{R}=\mathrm{Ind}^{\a_1, \ldots, \a_r}
\circ \mathrm{Res}^{\a_1, \ldots, \a_r}(\mathbb{P}) \in
\mathcal{Q}^{\alpha}.$$
As in \cite{S1}, Lemma~7.1., we have
\begin{equation}\label{E:proofcan1}
supp(\mathbb{R}) \subset
\bigcup_{(\underline{\beta}) \preceq (\underline{\alpha})} \underline{HN}^{-1}
(\underline{\beta})
\end{equation}
 and if
 $j: \underline{HN}^{-1}(\underline{\alpha}) \hookrightarrow \QQ^\alpha$
 is the embedding then
\begin{equation}\label{E:proofcan2}
j^*(\mathbb{R}) \simeq j^*(\mathbb{P}).
\end{equation}
Let us first assume that $\mu(\a_2)\neq \infty$, so that for all $i$ we have $rank(\alpha_i) < rank(\alpha)$. 
By the induction hypothesis we may find
for any $m \in \Z$ an element $w \in \mathfrak{W} ^{\otimes r}$ such that
$$\mathb{b}_{\text{Res}^{\a_1, \ldots, \a_r}(\mathbb{P})} \equiv w\; (\text{mod}\; (\widehat{\mathfrak{U}}^+[\a_1] \otimes \cdots \otimes \widehat{\mathfrak{U}}^+[\a_r])_{<m}),$$
where $$(\widehat{\mathfrak{U}}^+[\a_1] \otimes \cdots \otimes \widehat{\mathfrak{U}}^+[\a_r])_{<m}=\sum_{i=1}^r \widehat{\mathfrak{U}}_{\E,e}^+[\a_1] \otimes \cdots \otimes \widehat{\mathfrak{U}}^+[\a_i]_{<m} \otimes \cdots \otimes  \widehat{\mathfrak{U}}_{\E,e}^+[\a_r].$$
Since the induction product is continuous, we can choose $m \ll 0$ so that $\mathb{b}_{\mathbb{R}} =\mathb{b}_{\text{Ind}\circ \text{Res}(\mathbb{P})} \equiv \text{Ind}(w)\; (\text{mod}\; \widehat{\mathfrak{U}}^+_{<n})$.
Next, as $\mathbb{P}$ is simple, there exists a semisimple complex $\mathbb{P}'$ with $supp(\mathbb{P}')\subset \bigcup_{\underline{\beta} \prec \underline{\a}}
\underline{HN}^{-1}(\underline{\beta})$ such that $\mathbb{P} \oplus \mathbb{P}'\simeq \mathbb{R}$. Using the induction hypothesis again there exists $w' \in \mathfrak{W}$ such that $\mathb{b}_{\mathbb{P}'} \equiv w'\;(\text{mod}\; \widehat{\mathfrak{U}}^+_{<n})$, and we finally obtain $\mathb{b}_{\mathbb{P}} = \mathb{b}_{\mathbb{R}}-\mathb{b}_{\mathbb{P}'} \equiv \text{Ind}(w)-w' \;(\text{mod}\; \widehat{\mathfrak{U}}^+_{<n})$ as desired.
This closes the induction step when $r>1$ and $\underline{\a} \neq (\alpha_1, \a_2)$ with $\mu(\a_2)=0$. 
If we are in this last case then we have $rank(\alpha_1)=rank(\alpha)$ and we need a slightly different argument. Let us write
$\mathrm{Res}^{\a_1, \a_2}(\mathbb{P})=\sum_i \mathbb{T}_i \boxtimes \mathcal{Q}_i$. According to the first part of the proof,  $\mathb{b}_{\mathcal{Q}_i} \in \mathfrak{W}$ for any $i$.  Moreover since $deg(\a_1) < deg(\a)$ , there exists by the induction hypothesis elements $w_i \in \mathfrak{W}$ such that
$\mathb{b}_{\mathbb{T}_i} \equiv w_i \;(\text{mod}\; \widehat{\mathfrak{U}}^+_{<n})$. As $\widehat{\mathfrak{U}}^+_{<n}$ is stable under right multiplication by $\mathfrak{U}^+_{\E,e}[\b]$ for any $\b$ satisfying $\mu(\beta)=\infty$, we deduce that
$$\mathb{b}_{\text{Ind}\circ \text{Res}(\mathbb{P})} \equiv \sum_i w_i \mathb{b}_{\mathcal{Q}_i}\;(\text{mod}\; \widehat{\mathfrak{U}}^+_{<n})$$
and we may conclude the proof as in the case $\mu(\a_2)\neq \infty$ above.

\vspace{.15in}

\paragraph{} It remains to consider the situation of a generically semistable $\mathbb{P}$.
Write $\a=l\delta_{\mu}$ for some $l \geq 1$ and $\mu \in \mathbb{Q}$.
By Proposition~\ref{P:descript-perv} we have $\mathbb{P}=\mathbf{IC}(U^{(l\delta_{\mu})},\sigma)$ for some irreducible representation $\sigma$ of $\mathfrak{S}_l$.
Denote by $u_{\a}: \QQ^{(\a)} \hookrightarrow \QQ^\a$ the open embedding. From the existence of the maps $\rho_{\mu,\infty}$  (see Section~3.5.1.) and (\ref{E:semismall})
one deduces that 
\begin{equation}\label{E:stablesemismall}
\begin{split}
u_{l\delta_{\mu}}^*{\mathrm{Ind}}^{n_1\delta_{\mu}, \cdots
,n_r\delta_{\mu}}(\mathbf{IC}&(\underline{U}^{n_1\delta_{\mu}},{\sigma_1})
\boxtimes \cdots \boxtimes
\mathbf{IC}(\underline{U}^{n_r\delta_{\mu}},{\sigma_r}))\\
&=u_{l\delta_{\mu}}^*\mathbf{IC}(\underline{U}^{l\delta_{\mu}},ind_{\mathfrak{S}_{n_1}\times
\cdots \times
\mathfrak{S}_{n_r}}^{\mathfrak{S}_l}(\sigma_1 \otimes \cdots \otimes
\sigma_r)).
\end{split}
\end{equation}
On the other hand the element $\mathb{b}_{l\delta_{\mu}}$ belongs to $\mathfrak{W}$ for all $l \geq 1$.
As the Grothendieck group of $\mathfrak{S}_l$ is linearly spanned by the class of the trivial representation together with the classes ${ind}^{\mathfrak{S}_l}_{\mathfrak{S}_{n_1} \times \cdots \times \mathfrak{S}_{n_r}}(\sigma_1 \otimes \cdots \otimes \sigma_r)$ with $n_i >1$ for all $i$,
we conclude, using the induction hypothesis and the continuity of the product that there exists a
semisimple complex $\mathbb{T}$ which is \textit{not} generically semistable and 
such that $\mathb{b}_{\mathbb{T} \oplus \mathbb{P}} \in \mathfrak{W} + \widehat{\mathfrak{U}}^+_{<n}$. But then $supp(\mathbb{T}) \subset  \bigcup_{\underline{\beta} \prec \underline{\a}}
\underline{HN}^{-1}(\underline{\beta})$ and by the induction hypothesis again we have $\mathb{b}_{\mathbb{T}} \in \mathfrak{W} + \widehat{\mathfrak{U}}^+_{<n}$. Thus in the end $\mathb{b}_{\mathbb{P}} \in \mathfrak{W} + \widehat{\mathfrak{U}}^+_{<n}$ and statement $a)$ is proved for $\mathbb{P}$. This closes the induction argument and concludes the proof of Theorem~\ref{T:cano}. \qed  

\vspace{.2in}

\section{Geometric construction of $\mathbf{B}$}

\vspace{.1in}

In this section, we draw some of the consequences that Theorem~\ref{T:cano} and Corollary~\ref{C:cano} have regarding the canonical basis ${\mathbf{B}}$ of $\widehat{\UU}^+_{\Rb}$ considered in Section~2. 

\vspace{.2in}

\paragraph{\textbf{6.1.}} We will need one technical result. By Theorem~\ref{P:purity}, any $\mathbb{P} \in \mathcal{P}$ is endowed with a distinguished isomorphism $h_{\mathbb{P}}: \mathbb{P} \stackrel{\sim}{\to} \tilde{F}^*\mathbb{P}$ with respect to which it is pointwise pure of weight zero. Recall from Section~4 that we have set $v=-{q^{-1/2}}$ and more generally $v_f=-q^{-f/2}$ and let us denote by $\alpha,\overline{\alpha}$ the Frobenius eigenvalues in $H^1(\Eb, \qlb)$ \footnote{In \cite{BS}, the Frobenius eigenvalues are denoted $\sigma, \overline{\sigma}$. In this paper we reserve the latter notation for the formal parameters in the ring $\mathbf{R}=[\sigma^{\pm 1/2}, \overline{\sigma}^{\pm 1/2}]$.} . We also put $\tau=\alpha q^{-1/2}$ and $\tau_f=\tau^f$. Note that $|\tau|=1$.

\vspace{.1in}

\begin{prop}\label{P:verypure} Assume that $\tau\neq 1$.
Let $\mathbb{P} \in \mathcal{P}^\a$ and let $\mathcal{F} \in Coh(\E/\mathbb{F}_{q^e})$. The eigenvalues of $\tilde{F}^e$ acting on $H^i(\mathbb{P})_{|\mathbf{O}_{\mathcal{F}}}$ all belong to $q^{ei/2}\tau^{e\Z}$.\end{prop}
\noindent
\textit{Proof.} We will prove this for $e=1$; the other cases can be deduced by base change.
Fix $n \ll 0$ and let $x \in \mathbf{O}_{\mathcal{F},n}$. Let $\{\lambda_{ij}\}_{j \in J_i}$ be the set of Frobenius eigenvalues (with multiplicity) of $H^i(\mathbb{P}_n)_{|x}$. By Theorem~\ref{P:purity}, we may write $\lambda_{ij}=q^{i/2}\gamma_{ij} $ for some $\gamma_{ij}$ satisfying $|(\gamma_{i,j})|=1$, and by definition the coefficient of $[\mathcal{F}]$ in $Tr_f(\mathb{b}_{\mathbb{P}})$ is equal to
$$Tr_f(\mathb{b}_{\mathbb{P}})_{|\mathbf{O}_\mathcal{F}}=\sum_{i,j}(-1)^i\lambda_{ij}^f=\sum_{i,j}(-1)^{i(f+1)}v^{-if}\gamma_{ij}^f.$$
Our approach is based on the following elementary lemma~:

\begin{lem} Let $a, b$ be two complex numbers satisfying $|a|<1$, $|b|=1$ and $b \neq 1$.
Let $P(x,y) \in \mathbb{Q}[x^{\pm 1}, y^{\pm 1}]$ and assume given finite sets $K_i$ and complex numbers $\omega_{ij}, j \in K_i$ satisfying $|\omega_{ij}|=1$ such that, for any $f \geq 1$ 
$$P((-1)^{f+1}a^f,b^f)=\sum_i (-1)^{i(f+1)}a^{-if} \sum_{j \in K_i} \omega_{ij}^f.$$
Then $P(x,y) \in \N[x^{\pm 1},y^{\pm 1}]$ and $\omega_{ij} \in b^{\Z}$ for all $i,j$.\end{lem}
\noindent
\textit{Proof.} We argue by induction on the number $N(P)$ of monomials $x^iy^j$ appearing in $P(x,y)$.
If $N=0$, i.e. $P=0$ then $\sum_i \sum_{j \in K_i} (-1)^{i(f+1)}a^{-if}\omega_{ij}^f=0$ for all $f$. We will prove that this forces all the sets $K_i$ to be empty. Choose $i_0$ maximal such that $K_{i_0} \neq \emptyset$. Then 
\begin{equation}\label{E:sublem}
\begin{split}
0&=\sum_{i \leq i_0} (-1)^{i(f+1)}a^{(i_0-i)f} \sum_{j \in K_i}\omega_{ij}^f\\
&=(-1)^{i_0(f+1)}\sum_{j \in K_{i_0}} \omega_{i_0j}^f + \sum_{i < i_0} (-1)^{i(f+1)}a^{(i_0-i)f} \sum_{j \in K_i}\omega_{ij}^f
\end{split}
\end{equation}
for all $f$. Letting $f \mapsto \infty$ we get $\text{Lim}_{f \mapsto \infty} \sum_{j \in K_{i_0}} \omega_{i_0j}^f=0$. But this contradicts the following well known fact

\vspace{.05in}

\begin{sublem} Let $\{\omega_j\}$, $j \in K$ be a finite set of complex numbers satisfying $|\omega_j|=1$ for all $j$. If $K \neq \emptyset$ then $\sum_j \omega_j^N$ does not converge to zero as $N$ tends to infinity.
\end{sublem}

Now fix $P \neq 0$ and assume that the Lemma is proved for all $P'$ such that $N(P')<N(P)$. Let $x^{i_0}y^{j_0}$ be the monomial in $P$ of least degree, ordered lexicographically. Multiplying throughout by $(-1)^{(f+1)i_0}a^{-fi_0}b^{-fj_0}$ we may assume that $i_0=j_0=0$, i.e $P(x,y)=P_0(y) + \sum_{i>0} x^i P_i(y)$ and $P_0(y)=c_0 + \sum_{j>0}c_jy^j$. For any $l>0$ set
\begin{equation*}
\begin{split}
H_l&=\frac{1}{l} \sum_{f=1}^l P((-1)^{(f+1)}a^f, b^f)\\
&=c_0 + \sum_{j>0} c_j\big(\frac{1}{l}\sum_{f=1}^l b^{fj}\big)+\sum_{i>0} \frac{1}{l} \sum_{f=1}^l (-1)^{(f+1)i} a^{if}P_i(b^f).
\end{split}
\end{equation*}
As $|a|<1$ and $b \neq 1$ we have $\text{Lim}_{l \to \infty} H_l=c_0$. On the other hand $H_l =\sum_i  \sum_{j \in K_i} \frac{1}{l} \sum_{f=1}^l (-1)^{(f+1)i}a^{-if}\omega_{ij}^f$. This converges if and only if $K_i=\emptyset$ for $i>0$ and, noting that
$$\underset{\;l \to \infty}{\text{Lim}}\; \frac{1}{l} \sum_{f=1}^l \omega_{ij}^f=\begin{cases} 1 & \text{if}\; \omega_{ij}=1\\
0 & \text{otherwise}
\end{cases}$$
we have $c_0= \text{Lim}_{l \to \infty} H_l=\#\{j \in K_0\;|\; \omega_{0j}=1\}$. In particular, $c_0 \in \N$. We may now use the induction hypothesis with $P'(x,y)=P(x,y)-c_0$ and
$$K'_i=\begin{cases} K_i & \text{if}\; i \neq 0\\ K_0 \setminus \{j \in K_0\;|\; \omega_{0j}=1\}& \text{if}\; i=0
\end{cases}$$
to conclude that $P'(x,y)$, hence also $P(x,y)$ belongs to $\N[x^{\pm 1}, y^{\pm 1}]$. The proof also gives that $\omega_{ij} \in b^{\Z}$ for all $i,j$. We are done. \qed

\vspace{.1in}

\noindent
\textit{End of Proof of Proposition~\ref{P:verypure}.} By the above lemma, if a semisimple complex $\mathbb{R} \in \mathcal{Q}^{\a}$ satisfies the property

\vspace{.05in}

\noindent
(a) \textit{for all $\mathcal{F} \in Coh(\E)$ there exists $R\in \mathbb{Q}[x^{\pm 1},y^{\pm 1}]$ such that
$Tr_f(\mathb{b}_{\mathbb{R}})_{|\mathbf{O}_{\mathcal{F}}}=R(v_f,\tau_f)$ for all $f \geq 1$},

\vspace{.05in}

\noindent
then the Frobenius eigenvalues on all stalks of $H^i(\mathbb{R})$ belong to $q^{i/2}\tau^{\Z}$. The converse is also obviously true. We claim that if $\mathbb{R}$ and $\mathbb{R}'$ both satisfy (a) then so does $\mathrm{Ind}(\mathbb{R} \boxtimes \mathbb{R}')$. Indeed, if $\mathbb{R}$ and $\mathbb{R}'$ satisfy (a) then there exists
$R_i, R'_j \in \mathbb{Q}[x^{\pm 1}, y^{\pm 1}]$ and paths $\p_i, \p'_j \in \textbf{Conv}^+$ such that \begin{equation}\label{E:Trace}
Tr_f(\mathb{b}_{\mathbb{R}})=\sum_i R_i(v_f,\tau_f)\mathbf{1}^{\textbf{ss}}_{\p_i} , \qquad 
Tr_f(\mathb{b}_{\mathbb{R}'})=\sum_j R'_j(v_f,\tau_f)\mathbf{1}^{\textbf{ss}}_{\p'_j}.
\end{equation} 
Recall from \cite{BS}, Prop.~6.3, that $\widehat{\mathbf{U}}^+_{\E}$ is isomorphic to the specialization at $\sigma^{1/2}=\a^{1/2}, \overline{\sigma}^{1/2}=\overline{\alpha}^{1/2}$ of $\widehat{\boldsymbol{\mathcal{E}}}^+_{\mathbf{R}}$. Of course this holds for all the base field extensions $X/\mathbb{F}_{q^f}$ as well, with the corresponding specialization $\sigma^{1/2}=\a^{f/2}, \overline{\sigma}^{1/2}=\overline{\alpha}^{f/2}$. Set $\nu=-(\sigma \overline{\sigma})^{-1/2}, t=\sigma^{1/2}\overline{\sigma}^{-1/2}$. These specialize respectively to $v_f$ and $\tau_f$. We deduce from (\ref{E:Trace}) that 
$Tr_f(\mathb{b}_{\mathbb{R}})$ and $Tr_f(\mathb{b}_{\mathbb{R}'})$ are the specialization at $\nu=v_f$ and $t=\tau_f$ of ${u}=\sum_i R_i(\nu,t) {\mathbf{1}}^{\textbf{ss}}_{\p_i}$ and
${u}'=\sum_j R'_j(\nu,t) {\mathbf{1}}^{\textbf{ss}}_{\p'_j}$ respectively. But then
$Tr_f(\mathb{b}_{\mathrm{Ind}(\mathbb{R} \boxtimes \mathbb{R}')})$ is the specialization of
${u}\cdot {u}'$, and it follows that (a) holds for $\mathrm{Ind}(\mathbb{R} \boxtimes \mathbb{R}')$ as well.

\vspace{.05in}

\noindent
We may now finish the proof of Proposition~\ref{P:verypure}. It is clear that $\mathbbm{1}_{\a}$ satisfies (a) for any $\a$, and hence the same is true for products $\text{Ind}(\mathbbm{1}_{\a_1} \boxtimes \cdots \boxtimes \mathbbm{1}_{\a_r})$. We deduce that the Frobenius eigenvalues of any stalk of $H^i(\text{Ind}(\mathbbm{1}_{\a_1} \boxtimes \cdots \boxtimes \mathbbm{1}_{\a_r}))$ belong to
$q^{i/2}\tau^{\Z}$.
Observe that for any slope $\mu$ the simple perverse sheaf $\mathbf{IC}(U^{(l\delta_{\mu})},\sigma_{\lambda}) \in \mathcal{P}^{(l\delta_{\mu})}$ appears with multiplicity \textit{one} in $\text{Ind}(\mathbbm{1}_{\lambda_1\delta_{\mu}} \boxtimes \cdots \boxtimes \mathbbm{1}_{\lambda_r\delta_{\mu}})$ if $\lambda=(\lambda_1, \ldots, \lambda_r)$ and $\sigma_{\lambda}$ is the irreducible $\mathfrak{S}_l$-module of type $\lambda$. From this it follows that the Frobenius eigenvalues of any stalk of $H^i(\mathbf{IC}(U^{(l\delta_{\mu})},\sigma_{\lambda}))$
belong to $q^{i/2}\tau^{\Z}$. Thus $\mathbf{IC}(U^{(l\delta_{\mu})},\sigma_{\lambda})$ satisfies (a) in turn. Finally, by Proposition~\ref{P:descript-perv} any $\mathbb{P} \in \mathcal{P}$ appears with multiplicity one in some product $\text{Ind}(\mathbb{P}_1 \boxtimes \cdots \boxtimes \mathbb{P}_s)$ for certain $\mathbb{P}_i=\mathbf{IC}(U^{(l_i\delta_{\mu_i})},\sigma_{\lambda_i}))$. Arguing as above, we obtain that the Frobenius eigenvalues on $H^i(\mathbb{P})$ all belong to $q^{i/2}\tau^{\Z}$ as desired. We are done. \qed

\vspace{.2in}

\paragraph{\textbf{6.2.}} The proof of Proposition~\ref{P:verypure} gives us in fact the following property~:

\vspace{.05in}

 \textit{For any $\mathbb{P} \in \mathcal{P}$ there exists a unique ${\mathbf{b}}_{\mathbb{P}} \in \widehat{\boldsymbol{\mathcal{E}}}^+_{\Rb}$ such that for any $f \geq 1$,} 
 \begin{equation}\label{E:star}
 Tr_f(\mathb{b}_{\mathbb{P}})=(\mathbf{b}_{\mathbb{P}})_{\big|\substack{\sigma^{1/2}=\a^{f/2} \\ \overline{\sigma}^{1/2}=\overline{\alpha}^{f/2}}}.
 \end{equation}

\vspace{.1in}

By Proposition~\ref{P:descript-perv}, the set $\mathcal{P}^{(\a)}$ is indexed  by partitions of $deg(\a)$, and for $\mu(\a_1) < \cdots < \mu(\a_r)$, the correspondences
\begin{align*}
\mathcal{P}^{(\a_1)} \times \cdots \times \mathcal{P}^{(\a_r)} &\to \mathbf{Conv}^+\\
(\lambda_1^{(1)}, \ldots, \lambda_{s_1}^{(1)}; \cdots ; \lambda_1^{(r)}, \ldots, \lambda_{s_r}^{(r)}) &\mapsto \bigg(\lambda_1^{(1)}\frac{\a_1}{deg(\a_1)}, \ldots, \lambda_{s_r}^{(r)}\frac{\a_r}{deg(\a_r)}\bigg)
\end{align*}
all together set up a bijection $\omega: \mathcal{P} \stackrel{\sim}{\to} \mathbf{Conv}^+$.

\begin{lem}\label{L:Basis} The set $\{{\mathbf{b}}_{\mathbb{P}}\}$ is an $\Rb$-basis of
$\widehat{\boldsymbol{\mathcal{E}}}^+_{\Rb}$.\end{lem}
\noindent
\textit{Proof.} By construction we have, for $\x \in \ZZ^+$,
${\mathbf{b}}_{\mathbbm{1}_{\x}}={\mathbf{1}}_{\x}\in \beta_{\x} + \prod_{\p \prec \x} \Rb \beta_{\p}$. Moreover if $deg(\x)=1$ and $\sigma_i$ are representations of the symmetric group $\mathfrak{S}_{l_i}$ then
\begin{equation*}
\begin{split}
\mathrm{Ind}^{l_1\x, \ldots, l_r\x}\bigg(\mathbf{IC}(U^{(l_1\x)},\sigma_1)& \boxtimes \cdots
\boxtimes \mathbf{IC}(U^{(l_r\x)},\sigma_r)\bigg)\\
=&\mathbf{IC}\bigg(U^{(l\x)},ind_{\mathfrak{S}_{l_1} \times \cdots \times \mathfrak{S}_{l_r}}^{\mathfrak{S}_l} (\sigma_1 \times \cdots \times \sigma_r)\bigg) \oplus \mathbb{R},
\end{split}
\end{equation*}
where $l=\sum l_i$ and $Supp(\mathbb{R}) \subset \QQ^{l\x} \setminus \QQ^{(l\x)}$.
From the above we deduce that for such an $\x$ and for any partition 
$\lambda=(\lambda_1,\lambda_2,\ldots)$ of size $l$ we have 
${\mathbf{b}}_{\mathbb{P}}\in \beta_{\mathbf{p}} \oplus \prod_{\mathbf{q} \prec (l\x)} \Rb \beta_{\mathbf{q}}$
where $\mathbb{P}=\mathbf{IC}(U^{(l\x)}, \sigma_{\lambda}) \in \mathcal{P}^{(l\x)}$, $\sigma_{\lambda}$ is the irreducible $\mathfrak{S}_l$-module of type $\lambda$, and $\mathbf{p}=\omega(\mathbb{P})=(\lambda_1\x, \lambda_2\x,\ldots)$.
Finally, using this and Proposition~\ref{P:descript-perv} we obtain in turn that 
${\mathbf{b}}_{\mathbb{P}}\in \beta_{\omega(\mathbb{P})} \oplus \prod_{\mathbf{q} \prec \omega(\mathbb{P})} \Rb \beta_{\mathbf{q}}$ for an arbitrary $\mathbb{P}$.  As a consequence, $\{{\mathbf{b}}_{\mathbb{P}}\}$ is an $\Rb$-basis of $\widehat{\boldsymbol{\mathcal{E}}}^+_{\Rb}$ as wanted. \qed

\vspace{.1in}

\begin{cor} \label{C:verypure2} For any three simple complexes $\mathbb{P},\mathbb{P}',\mathbb{P}'' \in \mathcal{P}$ the eigenvalues of $\tilde{F}^e$ on $H^i(\text{Hom}(\mathbb{P}, \text{Ind}(\mathbb{P}' \boxtimes \mathbb{P}'')))$ all belong to $q^{ei/2}\tau^{e\Z}$. The same statement holds for
 $H^i(\text{Hom}(\mathbb{P}' \boxtimes \mathbb{P}'', \text{Res}(\mathbb{P}')))$.
\end{cor}
\noindent
\textit{Proof.} Let ${\mathbf{b}}_{\mathbb{P}}, {\mathbf{b}}_{\mathbb{P}'}, {\mathbf{b}}_{\mathbb{P}''}$ be the elements of $\widehat{\boldsymbol{\mathcal{E}}}^+_{\Rb}$ associated to $\mathbb{P},\mathbb{P}',\mathbb{P}''$. The coefficient of ${\mathbf{b}}_{\mathbb{P}}$ in the product ${\mathbf{b}}_{\mathbb{P}'}{\mathbf{b}}_{\mathbb{P}''}$ is given by a certain polynomial $K \in \mathbb{Q}[\nu^{\pm 1},t^{\pm 1}]$ and we have, for any $e \geq 1$,
$tr_e(\text{Hom}(\mathbb{P}, \text{Ind}(\mathbb{P}' \boxtimes \mathbb{P}'')))=K(v_e,\tau_e)$. The Corollary can now be deduced by the same argument as in Proposition~\ref{P:verypure}.\qed

\vspace{.2in}

\paragraph{\textbf{6.3.}} We presently come to the main result of this Section.

\begin{prop}\label{C:penultim} For any $\mathbb{P} \in \mathcal{P}$ we have ${\mathbf{b}}_{\mathbb{P}}={\mathbf{b}}_{\omega(\mathbb{P})}$. In particular,
\begin{enumerate}
\item[i)] For any $\mathbf{p} \in \mathbf{Conv}^+$ we have ${\mathbf{b}}_{\p} \in \beta_{\p} + \prod_{\p' \prec \p}\nu \mathbb{N}[\nu,t^{\pm 1}] \beta_{\p'}$,
\item[ii)] For any $\p,\p' \in \mathbf{Conv}^+$ we have
${\mathbf{b}}_{\p}{\mathbf{b}}_{\p'} \in \prod_{\p''} \mathbb{N}
[\nu^{\pm 1}, t^{\pm 1}] {\mathbf{b}}_{\p''}$.
\end{enumerate}
\end{prop}
\noindent
\textit{Proof.} Lemma~\ref{L:Basis} allows us to define in a unique way an involution $x \mapsto x^{\vee}$ of $\widehat{\boldsymbol{\mathcal{E}}}^+_{\Rb}$ by $\sigma^\vee=\sigma^{-1}, \overline{\sigma}^\vee=\overline{\sigma}^{-1}$  and ${\mathbf{b}}_{\mathbb{P}}^{\vee}={\mathbf{b}}_{\mathbb{P}}$. Note that $\nu^{\vee}=\nu^{-1}$, $t^{\vee}=t^{-1}$. We claim that this involution coincides with the bar involution $x \mapsto \overline{x}$ considered in Section~2. To see this, first observe that we may naturally extend the definition of ${\mathbf{b}}_{\mathbb{P}}$ to an arbitrary semisimple complex $\mathbb{P} \in \bigsqcup \mathcal{Q}^{\a}$ by bilinearity and by setting ${\mathbf{b}}_{\mathbb{P}[i](i/2)}=\nu^i{\mathbf{b}}_{\mathbb{P}}$. With this convention, (\ref{E:star}) holds for an arbitrary semisimple $\mathbb{P}$. Next, it is clear that ${\mathbf{b}}_{\mathrm{Ind}(\mathbb{P} \boxtimes \mathbb{P}')}={\mathbf{b}}_{\mathbb{P}}{\mathbf{b}}_{\mathbb{P}'}$ and  ${\mathbf{b}}_{\mathbb{P}}^{\vee}={\mathbf{b}}_{D(\mathbb{P})}$ (recall that by Proposition~\ref{P:descript-perv} every simple $\mathbb{P} \in \mathcal{P}$ is self-dual). 
It follows that $x \mapsto {x}^{\vee}$ is a ring homomorphism, and
in particular, we have ${\mathbf{1}}_{\a_1, \ldots, \a_r}^{\vee}={\mathbf{b}}_{\mathbbm{1}_{\a_1, \ldots, \a_r}}^{\vee}={\mathbf{b}}_{\mathbbm{1}_{\a_1}}^{\vee} \cdots {\mathbf{b}}_{\mathbbm{1}_{\a_r}}^{\vee}={\mathbf{b}}_{\mathbbm{1}_{\a_1}}\cdots {\mathbf{b}}_{\mathbbm{1}_{\a_r}}={\mathbf{1}}_{\a_1, \ldots, \a_r}$. We conclude that
$x \mapsto \overline{x}$ and $x \mapsto x^{\vee}$ coincide on the basis $\{{\mathbf{1}}_{\a_1, \ldots, \a_r}\}$, hence the two involutions are equal.

\vspace{.05in}

To obtain the equality ${\mathbf{b}}_{\mathbb{P}}={\mathbf{b}}_{\omega(\mathbb{P})}$ it remains to prove that 
\begin{equation}\label{E:vbase}
{\mathbf{b}}_{\mathbb{P}} \in \beta_{\omega(\mathbb{P})} + \prod_{\mathbf{q} \prec \omega(\mathbb{P})} \nu \mathbb{Q}[\nu,t^{\pm 1}] \beta_{\mathbf{q}}.
\end{equation}
Fix $\mathbb{P} \in \mathcal{P}^{\a}$. In order to show (\ref{E:vbase}) for $\mathbb{P}$, we will study the value of $Tr_e^{(n)}(\mathb{b}_{\mathbb{P}})$ over all points $x$ in $Q_n^{\a}$. Consider a HN type $\underline{\gamma}=(\gamma_1, \ldots, \gamma_r)$ of weight $\a$ and choose $x=(\phi: \mathcal{L}_n^{\a} \to \mathcal{F}) \in Q_n^{\a}$ with $HN(\mathcal{F})=\underline{\gamma}$. As in \cite{S1}, Lemma~7.1. we have
$$\mathbb{P}_{|\underline{HN}^{-1}(\underline{\gamma})}\simeq \text{Ind}^{\gamma_1, \ldots, \gamma_r} \text{Res}^{\gamma_1, \ldots, \gamma_r}(\mathbb{P})_{|\underline{HN}^{-1}(\underline{\gamma})}.$$

Let us write $\text{Res}^{\gamma_1, \ldots, \gamma_r}(\mathbb{P})=\bigoplus_i V_i \otimes \big( \mathbb{R}_1^i \boxtimes \cdots \boxtimes \mathbb{R}^i_r\big)$, where $\mathbb{R}_j^i \in \mathcal{P}^{\gamma_j}$ and $V_i$ is a complex of vector spaces. As $\mu(\gamma_1) < \cdots < \mu(\gamma_r)$ we deduce that, up to some (explicit) power of $v$,
$$Tr_1^{(n)}(\mathb{b}_{\mathbb{P}})(x)=\sum_i tr_1(V_i) Tr_1^{(n)}(\mathb{b}_{\mathbb{R}_1^i})(x_1) \times \cdots \times Tr_1^{(n)}(\mathb{b}_{\mathbb{R}_r^i})(x_r)$$
for any collection of semistable points $x_j=(\phi_j: \mathcal{L}_n^{\beta_j} \tto \mathcal{G}_j) \in Q_n^{(\gamma_j)}$ such that $\mathcal{F}\simeq \bigoplus_j \mathcal{G}_j$. Recall that by construction
$Tr_1^{(n)}(\mathb{b}_{\mathbb{R}_j^i})(x_j)=0$ unless $\mathbb{R}_j^i \in \mathcal{P}^{(\gamma_j)}$ is generically semistable, in which case $Tr_1^{(n)}(\mathb{b}_{\mathbb{R}_j^i})(x_j)={\beta}_{\omega(\mathbb{R}_j^i)}(x_j)$.
Using Proposition~\ref{P:verypure} and Corollary~\ref{C:verypure2} we deduce that
$$Tr_1^{(n)}(\mathb{b}_{\mathbb{P}})(x)\in\big(\beta_{\omega(\mathbb{P})} + \bigoplus_{\mathbf{q} \prec \omega(\mathbb{P})} \mathbb{N}[v^{\pm 1}, \tau^{\pm 1}] \beta_{\mathbf{q}}\big)(x).$$
A similar result holds if we replace $1$ by any $e \geq 1$ and $(v,\tau)$ by $(v_e,\tau_e)$. This implies that ${\mathbf{b}}_{\mathbb{P}}(x) \in \big(\beta_{\omega(\mathbb{P})} + \bigoplus_{\mathbf{q} \prec \omega(\mathbb{P})} \mathbb{N}[\nu^{\pm 1}, t^{\pm 1}] \beta_{\mathbf{q}}\big)(x)$. Since this is true for all $x \in Q_n^{\a}$ and all $n$, we finally obtain that ${\mathbf{b}}_{\mathbb{P}} \in \beta_{\omega(\mathbb{P})} + \bigoplus_{\mathbf{q} \prec \omega(\mathbb{P})} \mathbb{N}[\nu^{\pm 1}, t^{\pm 1}] \beta_{\mathbf{q}}$.

\vspace{.05in}

Finaly, recall that by Proposition~\ref{P:descript-perv} we have $\mathbb{P}=\mathbf{IC}(Z_{\mathbb{P}}, \mathfrak{L}_{\mathbb{P}})$ where if $\mathbb{P} \in \mathcal{P}^{\a_1} \times \cdots \times \mathcal{P}^{\a_l}$ then
$Z_{\mathbb{P}}$ is the substack of $\underline{Coh}^{\a}$ classifying coherent sheaves $\mathcal{H}$ isomorphic to a direct sum $\bigoplus_{i=1}^l \bigoplus_{j=1}^{deg(\a_i)} \mathcal{H}_{ij}$ where $\mathcal{H}_{i,j}$ is a stable sheaf of class $ \a_i/deg(\a_i)$, and $\mathfrak{L}_{\mathbb{P}}$ is a certain local system on $Z_{\mathbb{P}}$. This substack is also locally of finite type, and is equal to the increasing union of quotient stacks $Z_{\mathbb{P},n}/G_n^{\a}$ where
\begin{equation}\label{E:XP}
\begin{split}
Z_{\mathbb{P},n}= \{\phi: \mathcal{L}_n^{\a} \tto \mathcal{H}\;|\; \mathcal{H} \simeq \bigoplus_{i=1}^l \bigoplus_{j=1}^{deg(\a_i)} \mathcal{H}_{ij},\;\text{with}\;
 \mathcal{H}_{ij}\;\text{a\;stable\;sheaf\;of\;class\;} \a_i/deg(\a_i)\}
\end{split}
\end{equation}
and $\mathbb{P}_n=\mathbf{IC}(Z_{\mathbb{P},n}, \mathfrak{L}_{\mathbb{P}})$.
By Proposition~\ref{P:verypure} there exists, for any fixed $y \in Q_n^{\a}$, Laurent polynomials $P_{i,y}(t) \in \mathbb{N}[t^{\pm 1}]$ such that 
$$Tr_1^{(n)}(\mathb{b}_{\mathbb{P}})(y)=v^{-\dim\;G_n^{\a}}\sum_i (-1)^i tr_1(H^i(\mathbb{P}_n)_{|y})=v^{-\dim\;G_n^{\a}}\sum_i v^{-i}P_{i,y}(t).$$
By definition of the intersection cohomology complex, we have
$$(\mathbb{P}_n)_{|Z_{\mathbb{P},n}}=\mathfrak{L}_{\mathbb{P}}[\dim\;Z_{\mathbb{P},n}], \qquad
\dim\;supp(H^{-i}(\mathbb{P}_n)) <i\qquad\;\text{for\;} i < \dim\; Z_{\mathbb{P},n}.$$
In particular, as $\mathbb{P}_n$ is locally constant on each of the subvarieties $Z_{\mathbb{Q},n}$ we have
$P_{j,y}=0$ for any pair $(j,y)$ such that $y \in Z_{\mathbb{Q},n}$ with $\mathbb{Q} \neq \mathbb{P}$ and $j \geq -\dim\;Z_{\mathbb{Q},n}$. Hence the value of $Tr_1^{(n)}(\mathb{b}_{\mathbb{P}})$ at any $y \in Z_{\mathbb{Q},n}$ with $\mathbb{Q} \neq \mathbb{P}$ is the specialization at $\nu=v, t=\tau$ of an element in $\nu^{\dim\;Z_{\mathbb{Q},n}+1-\dim\;G_n^{\a}} \mathbb{N}[\nu,t^{\pm 1}]$. On the other hand, one checks that for any such $y$
$$\beta_{\omega(\mathbb{Q})}(y) \in \nu^{\dim\;Z_{\mathbb{Q},n}-\dim\;G_n^{\a}}\mathbb{N}[\nu], \qquad
\beta_{\omega(\mathbb{R})}(y)=0\;\text{if}\; \mathbb{R} \neq \mathbb{Q}.$$
We conclude that $\beta_{\omega(\mathbb{Q})}$ appears in ${\mathbf{b}}_{\mathbb{P}}$ with a coefficient in $\nu\mathbb{N}[\nu,t^{\pm 1}]$. Since this is true for any $\mathbb{Q}$, (\ref{E:vbase}) holds, and ${\mathbf{b}}_{\mathbb{P}}={\mathbf{b}}_{\omega(\mathbb{P})}$. Finally, statement i) was shown in the course of the above proof, and ii) is a consequence of Corollary~\ref{C:verypure2}.
\qed

\vspace{.1in}

\noindent
\textbf{Remark.} From the above proof it follows that the involution $x \mapsto \overline{x}$ of $\widehat{\mathbf{A}}^+_{\mathcal{A}}$ considered in Section~2 is a ring homomorphism.

\vspace{.2in}

\paragraph{\textbf{6.4.}}
It seems natural at this point to introduce a famlily of polynomials $\tilde{\daleth}_{\mathbf{p},\mathbf{q}}(\nu,t) \in \mathbb{N}[\nu,t^{\pm 1}]$ indexed by pairs of convex paths by the formula 
${\mathbf{b}}_{\mathbf{p}}=\sum_{\mathbf{q}} \tilde{\daleth}_{\mathbf{p},\mathbf{q}}(\nu,t) \beta_{\mathbf{q}}$. By construction we have $\tilde{\daleth}_{\mathbf{p},\mathbf{p}}(\nu,t)=1$ and $\tilde{\daleth}_{\mathbf{p},\mathbf{q}}(\nu,t)=0$ if $\mathbf{q} \not\prec \mathbf{p}$ and $\mathbf{q}\neq \mathbf{p}$.
In order to get something more reminiscent of Kostka polynomials, we slightly renormalize these polynomials as follows. Define a new basis $(\rho_{\mathbf{p}})_{\mathbf{p} \in \mathbf{Conv}^+}$ of $\widehat{\boldsymbol{\mathcal{E}}}^+_{\Rb}$ in exactly the same way as $(\beta_{\mathbf{p}})$ was defined, but using the Hall-Littlewood polynomials $P_{\lambda}(\nu)$ instead of the Schur function $s_{\lambda}$ (see Section 2.3.). We define the \textit{elliptic Kostka polynomial} $\daleth_{\mathbf{p},\mathbf{q}}(\nu,t)$ by the relations ${\mathbf{b}}_{\mathbf{p}}=\sum_{\mathbf{q}} {\daleth}_{\mathbf{p},\mathbf{q}}(\nu,t) \rho_{\mathbf{q}}$. Again, we have ${\daleth}_{\mathbf{p},\mathbf{p}}(\nu,t)=1$ and now ${\daleth}_{\mathbf{p},\mathbf{q}}(\nu,t)=0$ if $\mathbf{q} \not\preceq \mathbf{p}$. In addition, if
$\lambda=(\lambda_1, \ldots, \lambda_r)$ and $\sigma=(\sigma_1, \ldots, \sigma_s)$ are partitions and if $deg(\x)=1$ then setting
$\mathbf{p}=(\lambda_1\x, \ldots, \lambda_r\x), \mathbf{q}=(\sigma_1\x, \ldots, \sigma_s \x)$ we have $\daleth_{\mathbf{p},\mathbf{q}}(\nu,t)=
K_{\lambda,\sigma}(\nu)$ , where $K_{\lambda,\sigma}(\nu)$ denotes the usual q-Kostka polynomial.

\vspace{.1in}

We finish by describing a symmetry property which the polynomials $\tilde{\daleth}_{\mathbf{p},\mathbf{q}}(\nu,\tau)$ and $\daleth_{\mathbf{p},\mathbf{q}}(\nu,t)$ enjoy.

\begin{prop} Let $\gamma \in SL(2,\Z)$, and assume that $\mathbf{p},\mathbf{q} \in \mathbf{Conv}^+$ are such that $\gamma(\mathbf{p}), \gamma(\mathbf{q}) \in \mathbf{Conv}^+$ (that is, $\gamma(\mathbf{p})$ and $\gamma(\mathbf{q})$ are still entirely contained in $\ZZ^+$). Then
$$\daleth_{\gamma(\mathbf{p}),\gamma(\mathbf{q})}(\nu,t)=\daleth_{\mathbf{p},\mathbf{q}}(\nu,t), \qquad
\tilde{\daleth}_{\gamma(\mathbf{p}),\gamma(\mathbf{q})}(\nu,t)=\tilde{\daleth}_{\mathbf{p},\mathbf{q}}(\nu,t).$$
\end{prop}
\noindent
\textit{Proof.} By \cite{BS} Corollary~3.2 the group $SL(2,\Z)$ naturally acts on the Drinfeld double (with trivial center) $\boldsymbol{\mathcal{E}}_{\Rb}$ of $\boldsymbol{\mathcal{E}}_{\Rb}$. This action is compatible with the elements ${\mathbf{1}}^{\ss}_{\mathbf{p}}$ and $\beta_{\mathbf{p}}$; if $\gamma \in SL(2,\Z)$ is such that $\gamma(\mathbf{p}) \in \mathbf{Conv}^+$ then $\gamma \cdot {\mathbf{1}}^{\ss}_{\mathbf{p}}={\mathbf{1}}^{\ss}_{\gamma(\mathbf{p})}$ and $\gamma \cdot \beta_{\mathbf{p}}=\beta_{\gamma(\mathbf{p})}$. Now fix $\gamma,\mathbf{p},\mathbf{q}$ as in the 
Proposition. Recall from Section~2.2 that for a convex path $\mathbf{r}$ for which $\omega^{-1}(\mathbf{r}) \in \mathcal{P}^{\a_1} \times \cdots \times \mathcal{P}^{\a_r}$ we set $HN(\mathbf{r})=(\a_1, \ldots, \a_r)$ and that $\mathbf{r} \prec \mathbf{r}'$ is equivalent to $HN(\mathbf{r}) \prec HN(\mathbf{r}')$.
The action of $\gamma$ on $\boldsymbol{\mathcal{E}}_{\Rb}$ lifts at the geometric level to an isomorphism between open substacks of $\QQ^{|\mathbf{p}|}$ and $\QQ^{|\gamma(\mathbf{p})|}$~:
$$i_{\gamma}: \bigsqcup_{HN(\underline{\sigma}) \succ HN(\mathbf{q})} \underline{HN}^{-1}(\underline{\sigma}) \stackrel{\sim}{\longrightarrow} \bigsqcup_{HN(\underline{\sigma'}) \succ HN(\gamma(\mathbf{q}))} \underline{HN}^{-1}(\underline{\sigma})$$
which maps isomorphically the substacks $Z_{\omega^{-1}(\mathbf{p})}$ onto $Z_{\omega^{-1}(\gamma(\mathbf{p}))}$ (see \ref{E:XP}).
In particular, this allows us to identify the Frobenius traces of stalks of $\omega^{-1}(\mathbf{p})$ on
$Z_{\omega^{-1}(\mathbf{q})}$ and of $\omega^{-1}(\gamma(\mathbf{p}))$ on
$Z_{\omega^{-1}(\gamma(\mathbf{q}))}$. The Proposition follows.\qed

\vspace{.2in}

\centerline{\textbf{Acknowledgments}}

\vspace{.1in}

I would like to thank I. Burban, D. Harari, D. Madore and E. Vasserot for many helpful discussions.

\vspace{.2in}

\small{}

\end{document}